\documentclass{article}

\usepackage{amsmath,amsthm,amssymb}
\usepackage{amsmath}

\usepackage{stackrel}

\usepackage{cite}

\newtheorem{theorem}{Theorem}[section]

\newtheorem{corollary}[theorem]{Corollary}

\newcommand{\NN}{\mathbb{N}}
\newcommand{\ZZ}{\mathbb{Z}}

\newcommand{\RR}{\mathbb{R}}
\newcommand{\CC}{\mathbb{C}}

\newcommand{\ph}{\varphi}

\title{The concept of ``character'' in
Dirichlet's theorem\\ on primes in an arithmetic progression\footnote{This essay draws extensively on the second author's Carnegie Mellon MS thesis \cite{morris:11}. We are grateful to Michael Detlefsen and the participants in his \emph{Ideals of Proof} workshop, which provided feedback on portions of this material in July, 2011, and to Jeremy Gray for helpful comments. Avigad's work has been partially supported by National Science Foundation grant DMS-1068829 and Air Force Office of Scientific Research grant FA9550-12-1-0370.}}

\author{Jeremy Avigad and Rebecca Morris}

\begin{document}

\maketitle

\begin{abstract}
 In 1837, Dirichlet proved that there are infinitely many primes in any arithmetic progression in which the terms do not all share a common factor. We survey implicit and explicit uses of \emph{Dirichlet characters} in presentations of Dirichlet's proof in the nineteenth and early twentieth centuries, with an eye towards understanding some of the pragmatic pressures that shaped the evolution of modern mathematical method.
\end{abstract}

\tableofcontents

\section{Introduction}
\label{introduction:section}

Historians commonly take the ``modern'' age of mathematics to have begun in the nineteenth century. But although there is consensus that the events of that century had a transformational effect on mathematical thought, it is not easy to sum up exactly what changed, and why. Aspects of the transformation include an increasingly abstract view of mathematical objects; the rise of algebraic methods; the unification of disparate branches of the subject; evolving standards of rigor in argumentation; a newfound boldness in dealing with the infinite; emphasis on ``conceptual'' understanding, and a concomitant deemphasis of calculation; the use of (informal) set-theoretic language and methods; and concerns to identify a foundational basis to support the new developments.\footnote{The essays in Ferreir\'os and Gray \cite{ferreiros:gray:06} provide an overview.} It is still an important historical and philosophical task to better understand these components, and the complex interactions between them.

A good deal of attention has been given to early appearances of set-theoretic and structural language, including the use of equivalence relations and  ideals in algebra; the expansion of the function concept in analysis, and its generalization to other mathematical domains; foundational constructions of number systems from the natural numbers to the reals and beyond; and the overall conception of mathematics as the study of structures and spaces, often characterized in set-theoretic terms \cite{ferreiros:99}. These provide clear and focused manifestations of the changes that took place.

In this essay, we will consider certain functions, known as ``Dirichlet characters,'' and their role in proving a seminal 1837 theorem due to Dirichlet, which states that there are infinitely many prime numbers in any arithmetic progression in which the terms do not all share a common factor. As far as abstract objects go, characters are fairly benign: for a given positive integer $m$, there are only finitely many characters modulo $m$, and each one can be described exactly by giving its value on the finitely many residue classes modulo $m$. Moreover, it is not terribly hard to provide a concrete and exhaustive description of the set of all such characters. 

Nonetheless, we will argue that the evolution of the treatment of characters over the course of the nineteenth century illustrates important themes in the overall transformation of mathematical thought. Indeed, we will try to call attention to changes in mathematical method that have by now become so ingrained that it is hard for us today to appreciate their significance. An important feature of contemporary proofs of Dirichlet's theorem is that they are higher-order, which is to say, one treats characters on par with mathematical objects like the natural numbers: one sums over characters, quantifies over characters, applies functions to characters, forms sets of characters, considers groups whose elements are characters, and so on. This runs counter to an ``intensional'' view of characters as \emph{expressions} that play a role in ordinary mathematical language that is fundamentally different from that of the natural numbers or the real numbers. Our narrative traces a gradual transition from this latter view to the contemporary one.

In tracing this evolution, we will set aside other important aspects of the history of Dirichlet's theorem, such as the evolution of number-theoretic and analytic ideas, and the use of analysis in number theory. We will also for the most part set aside other nineteenth century developments that invoked the use of set-theoretic language. Our hope is that focusing on the notion of character, and its evolving role in a single mathematical proof, will help clarify some of the mathematical pressures that contributed to the adoption of a modern standpoint.

The outline of this paper is as follows. In Section~\ref{functions:section}, we situate our study amidst a host of topics related to the development of the function concept, noting that for most of the nineteenth century the word ``function'' was used exclusively in connection with functions on the real or complex numbers. In Section~\ref{contemporary:section}, we discuss contemporary presentations of Dirichlet's proof, and in Section~\ref{objects:section}, we highlight various senses in which these presentations treat characters as objects. In contrast, Section~\ref{dirichlet:section} describes Dirichlet's own presentation of his proof, in which the notion of a character does not figure at all. Section~\ref{evolution:section} then traces a gradual transition, as characters are transformed from shade-like entities in the original proof to the fully embodied objects we take them to be today. Section~\ref{analysis:section} analyzes the forces that shaped the transition, and Section~\ref{conclusions:section} draws some conclusions. In a companion paper to this one \cite{avigad:morris:b}, we explore further philosophical aspects of the treatment of characters as objects in the history of Dirichlet's theorem.

\section{Functions in the nineteenth century}
\label{functions:section}

\subsection{The generalization of the function concept}
\label{generalization:section}

Nineteenth century mathematicians dealt with a number of objects that we, today, view as instances of the function concept. Analysis dealt with functions defined on the continuum, that is, functions from the real numbers, $\RR$, to $\RR$. Cauchy and others extended the subject to include functions from the complex numbers, $\CC$, to $\CC$. Analysts also commonly dealt with sequences and series, which can be viewed as functions from the natural numbers, $\NN$, to a target domain, such as the integers or the reals. Number theory dealt with functions like the Euler totient function, $\ph$, which maps $\NN$ to $\NN$. Geometry dealt with transformations of the Euclidean or projective plane, which can be viewed as invertible maps from the space to itself. Galois' theory of equations focused on substitutions, or permutations, of a finite set $A$, which is to say, bijective functions from $A$ to $A$. Below we will see how developments in number theory led to the study of characters, that is, functions from the integers, $\ZZ$, to $\CC$, or from the group $(\ZZ/m\ZZ)^*$ of units modulo $m$ to $\CC$. By the end of the century, Dedekind and Cantor considered arbitrary mappings, or correspondences, between domains.

Despite this diversity, most of the literature on the evolution of the function concept (including \cite{chorlay:unp,Monna1972,Y1976,Kleiner1989,Luzin98}) has focused on functions on the real or complex numbers. Surveys typically trace the evolution of the concept from the introduction of the term ``function'' (``functio'') by Euler \cite{Euler48} in 1748, to the mention of ``arbitrary functions'' (``fonction arbitraire'') in the title of Dirichlet's seminal paper of 1829 \cite{Dirichlet1829}, to the dramatically opposed treatment of the function concept by Riemann and Weierstrass.\footnote{See Bottazzini and Gray \cite{bottazzini:gray:13} for a history of complex function theory, which includes a detailed treatment of the work of Riemann and Weierstrass, in particular.} There is a very good reason for this: in the early nineteenth century, the word ``function'' was used almost exclusively in this narrow sense. Indeed, there is little in the early nineteenth century literature that suggests any family resemblances between the kinds of mathematical objects enumerated in the last paragraph, let alone subsumption under a single overarching concept.

For example, sequences of real numbers were referred to as sequences or series, or often just introduced by displaying the initial elements $a_0, a_1, a_2, \ldots$ with an ellipsis at the end. The words ``series'' (``s\'erie'' in French, ``Reihe'' in German) was also used to describe a finite or infinite sum. For example, Dirichlet wrote the following in 1835:
\begin{quote}
Now let
\[
 F(\alpha) = b_0 + b_1 \cos \alpha + b_2 \cos 2\alpha + \ldots = \sum b_i \cos i \alpha
\]
be an arbitrary finite or infinite series, whose coefficients are independent of $\alpha$.\footnote{``Es sei jetzt: \[
 F(\alpha) = b_0 + b_1 \cos \alpha + b_2 \cos 2\alpha + \ldots = \sum b_i \cos i \alpha
\]
eine beliebige endliche oder unendliche Reihe, deren Coefficienten von $\alpha$ unabh\"{a}ngig sind.'' } \cite[p.~249]{dirichlet:35}
\end{quote}
The dependence of the coefficient $b_i$ on $i$ was indicated by the subscript, just as we do today.

Terminology governing number-theoretic functions was more varied. In the eighteenth century, Euler tended to refer to number-theoretic functions as ``symbols'' or ``characters.''  When introducing what we now call the totient function in 1781, he wrote:
\begin{quote}
 \ldots let the character $\pi D$ denote that multitude of numbers which are less than $D$ and which have no common divisor with it.\footnote{``Quod quo facilius praestari possit, denotet character $\ph D$ multitudinem istam numerorum ipso $D$ minorum, et qui cum eo nullum habeant divisorem communem.'' \cite[p.~19]{Euler564}} \cite[p.~2]{Euler564}
\end{quote}
In 1801, in \S 38 of the \emph{Disquisitiones Arithmeticae} \cite{gauss:01}, Gauss introduced the totient function as a notation:
\begin{quote}
 For brevity we will designate the number of positive numbers which are relatively prime to the given number and smaller than it by the prefix $\ph$. We seek therefore $\ph A$.\footnote{``Designemus brevitatis gratia multitudinem numerorum positivorum ad numerum datum primorum ipsoque minorum per praefixum characterem $\ph$.  Quaeritur itaque $\ph A$.''}
\end{quote}
In the following section, he referred to the ``character'' $\ph$. Later, in \S 52, he defined another number-theoretic function, $\psi$, and  in \S 53 compared it to the $\ph$ ``symbol'' (``signo''). In the \emph{Disquisitiones}, the word ``function'' is never used to describe such entities.


Early in the nineteenth century, in the field of projective geometry, M\"obius studied general ``relations'' between figures in the plane, such as ``affinities'' and ``collineations.''\footnote{See the discussion in Wussing  \cite[pp.~35--40]{wussing:84}.} In his \emph{Erlangen Program} of 1872, Klein considered invertible ``transformations'' (``Transformationen'') of a space, and the group they form under the operation of ``composition'' (``Zusammensetzung'') \cite{klein:93}. Once again, in neither case is there any hint that such transformations bear a connection to a more general function concept.

In algebra, Galois' groups consisted of ``substitutions'' of the roots. Cayley's famous 1854 paper, which provided the first axiomatization of the group concept, begins as follows:
\begin{quote}
 Let $\theta$ be a symbol or operation, which may, if we please, have for its operand, not a single quantity $x$, but a system $(x, y, \ldots)$, so that
\[
 \theta(x, y, \ldots) = (x', y', \ldots),
\]
 where $x', y', \ldots$ are any functions whatever of $x, y, \ldots$, it is not even necessary that $x', y', \ldots$ should be the same in number with $x, y, \ldots$. In particular, $x', y'$, \&c.~may represent a permutation of $x$, $y$, \&c., $\theta$ is in this case what is termed a substitution; and if, instead of a set $x, y, \ldots$, the operand is a single quantity $x$, so that $\theta x = x' = f x$, $\theta$ is an ordinary function symbol.\footnote{This excerpt is quoted and discussed in Pengelley~\cite{pengelley:05}.} \cite[p.~123]{cayley:54}
\end{quote}
Notice that although we can assume that each quantity $x', y', \ldots$ is a function of $x, y, \ldots$, for Cayley, $\theta$ is, in general, an operation and not a function. As examples of cases where $\theta$ \emph{is} a function, he singled out multiplication by quaternions and examples arising from the study of elliptic functions. But he makes it clear that it is the operations, more generally, that are assumed to satisfy the associative law, under composition:
\begin{quote}
 A symbol $\theta \ph$ denotes the compound operation, the performance of which is equivalent to the performance, first of the operation $\ph$, and then of the operation $\theta$; $\theta \ph$ is of course in general different from $\ph \theta$. But the symbols $\theta, \ph, \ldots$ are in general such that $\theta.\ph \chi = \theta \ph.\chi$, \&c., so that $\theta \ph \chi, \theta\ph\chi\omega$, \&c.~have a definite signification independent of the particular mode of compounding the symbols\ldots
\end{quote}

Dedekind's later treatment of Galois theory \cite{dedekind:94} focused on the group of automorphisms of a field rather than permutations of roots, but even in this context he used the term ``substitution'' for an invertible map from a field to itself, and ``permutation'' for a substitution that is moreover a homomorphism with respect to the field structure. Once again, there seems to be nothing to link these substitutions and permutations with the function concept from analysis. Even as late as 1895, when Cantor presented his mature theory of the infinite \cite{cantor:95}, he described a one-to-one correspondence as a ``law of association'' (``Zuordnungsgesetz'') between sets, with nothing to suggest that these correspondences had anything to do with the objects of function theory with which he had begun his mathematical career.

The shift from using the term ``function'' exclusively for functions defined on the continuum to those defined on more general domains was gradual. We have found a very early instance of the phrase ``number-theoretic function'' (``zahlentheoretische Funktion'') in an 1850 paper by Eisenstein, which begins with a notably self-conscious justification for the use of the term:
\begin{quote}
  Since, with the concept of a function, one moved away from the necessity of having an analytic construction, and began to take its essence to be a tabular collection of values associated to the values of one or several variables, it became possible to take the concept to include functions which, due to conditions of an arithmetic nature, have a determinate sense only when the variables occurring in them have integral values, or only for certain value-combinations arising from the natural number series.  For intermediate values, such functions remain indeterminate and arbitrary, or without any meaning.\footnote{``Seit man bei dem Begriffe der Funktion von der Nothwendigkeit der analytischen Zusammensetzung abgehend, das Wesen derselben in die tabellarische Zusammenstellung einer Reihe von zugeh\"{o}rigen Werthen mit den Werthen des order der (mehrerer) Variabeln zu setzen anfing, war es m\"{o}glich, auch solche Funktionen unter diesen Begriff mit aufzunehmen, welche aus Bedingungen arithmetischer Natur entspringend nur f\"{u}r ganze Werthe oder nur f\"{u}r gewisse aus der nat\"{u}rlichen Zahlenreihe hervorgehende Werthe und Werth-Combinationen der in ihnen vorkommenden Variabeln einen bestimmten Sinn erhalten, w\"{a}hrend sie f\"{u}r die Zwischenwerthe entweder unbestimmt und willk\"{u}rlich oder ohne alle Bedeutung bleiben.'' We are grateful to Wilfried Sieg for help with the translation.} \cite[p.~706]{eisenstein:50}
\end{quote}
We might have expected Eisenstein to observe, more simply, that once one starts to think of a function as a correspondence between input and output values, it is reasonable to transfer the notion to correspondences between domains other than the real or complex numbers, described by means other than an analytic expression. Instead, he adopted the surprising strategy of viewing number-theoretic functions as \emph{partial functions on the real numbers}. This is a nice illustration of the fact that concepts are often stretched and transformed in suprising ways, a theme that is the central study of Mark Wilson's book, \emph{Wandering Significance} \cite{wilson:06}. In any case, the passage makes it clear that the modern notion of a function defined on arbitrary domains was far from the mid-nineteenth century mindset.

We have already noted that in 1854, Cayley considered multiplication by a quaternion as a function. This seems to fit with the view of the natural progression from real numbers to complex numbers to quaternions, as generalizations of the concept of magnitude. One finds another hint of expansion in the title of Dedekind's 1854 \emph{Habilitation} lecture, ``On the introduction of new functions in mathematics''(``\"Uber die Einf\"uhrung neuer Funktionen in der Mathematik'') \cite{dedekind:54}. In this lecture, Dedekind discussed the way that the domain of natural numbers was gradually expanded to include the integers, real numbers, and complex numbers, while extending and preserving the properties of familiar operations like addition, and division. But the evocative use of the word ``function'' in the title is tempered by the contents of the lecture itself, where the word ``operation'' is used exclusively when it comes for functions defined on the natural numbers and integers.

To our knowledge, it is not until 1879 that one finds the word ``function'' used to describe correspondences between arbitrary domains; and then, in that year, this occurred in two remarkable sources. The first is due to Dedekind. In 1879, in a supplement dealing with quadratic forms in the third edition of his presentation of Dirichlet's lectures on number theory \cite{dirichlet:63b}, he defined the general notion of a character on a class group (that is, a function which maps equivalence classes of ideals in an algebraic number field to the complex numbers):\begin{quote}
 \ldots the function $\chi(\mathfrak a)$ also possesses the property that it takes the same value on all ideals $\mathfrak a$ belonging to the same class $A$; this value is therefore appropriately denoted by $\chi(A)$ and is clearly always an $h$th root of unity. Such functions $\chi$, which in an extended sense can be termed \emph{characters}, always exist; and indeed it follows easily from the theorems mentioned at the conclusion of \S 149 that the class number $h$ is also the number of all distinct characters $\chi_1, \chi_2, \ldots, \chi_h$ and that every class $A$ is completely characterized, i.e.~is distinguished from all other classes, by the $h$ values $\chi_1(A), \chi_2(A), \ldots, \chi_h(A)$.\footnote{``\ldots die Function $\chi(\mathfrak a)$ ausser der Eigenschaft \ldots noch die andere besitzt, f\"{u}r alle derselben Classe $A$ angeh\"{o}renden Ideale $\mathfrak a$ denselben Werth anzunehmen, welcher mithin zweckm\"{a}ssig durch $\chi(A)$ bezeichnet wird und offenbar immer eine $h^{te}$ Wurzel der Einheit ist.  Solche Functionen $\chi$, die man im erweiterten Sinn \emph{Charaktere} nennen kann, existiren immer, und zwar geht aus den am Schlusse des $\S 149$ erw\"{a}hnten S\"{a}tzen leicht hervor, dass die Classenanzahl $h$ zugleich die Anzahl aller verschiedenen Charaktere $\chi_{1}, \chi_{2}, \ldots, \chi_{h}$ ist, und dass jede Classe $A$ durch die ihr entsprechenden $h$ Werthe $\chi_{1}(A), \chi_{2}(A), \ldots, \chi_{h}(A)$ vollst\"{a}ndig charakterisirt, d.h.~von allen anderen Classen unterschieden wird.'' The quotation appears in \S 178 in the 1879 edition of the \emph{Vorlesungen} \cite{dirichlet:63b}, and in \S 184 of the 1894 edition, which is reproduced in Dedekind's \emph{Werke} \cite{dedekind:68}. The translation above is by Hawkins \cite[p.~149]{hawkins:71}.}
\end{quote}
In other words, Dedekind referred to a mapping from a certain algebraic structure to the complex numbers as a ``function.'' There is no reason to think that he meant to limit the terminology to this particular structure, and, indeed, by 1882, Weber \cite{weber:82} described characters on arbitrary groups as ``functions.''\footnote{We discuss Weber's 1882 paper, and provide more background on the history of characters, in Sections~\ref{dedekind:section} and \ref{weber:section} below.} In his foundational essay \emph{Was sind und was sollen die Zahlen} \cite{dedekind:88} of 1888, Dedekind considered arbitrary mappings (``Abbildung'') between domains, but, curiously, seems to distinguish mappings from functions.\footnote{See \S 135 in \cite{dedekind:88}. The relationship between the two notions will be discussed in forthcoming work by Wilfried Sieg and the second author.}

The other landmark source of 1879 is the philosopher Gottlob Frege's \emph{Begriffsshrift} \cite{Frege1879}.  Here Frege defined the notion of a function as follows:
\begin{quote}
If, in an expression (whose content need not be a judgeable content), a simple or complex symbol occurs in one or more places, and we think of it as replaceable at all or some of its occurrences by another symbol (but everywhere by the same symbol), then we call the part of the expression that on this occasion appears invariant the function, and the replaceable part its argument.\footnote{``Wenn in einem Ausdrucke, dessen Inhalt nicht beurtheilbar zu sein braucht, ein einfaches oder zusammengesetztes Zeichen an einer oder an mehren Stellen vorkommt, und wir denken es an allen oder einigen dieser Stellen durch Anderes, \"{u}berall aber durch Dasselbe ersetzbar, so nennen wir den hierbei unver\"{a}nderlich erscheinenden Theil des Ausdruckes Function, den ersetzbaren ihr Argument.''} \cite[\S 9]{Frege1879}
\end{quote}
Although this definition is couched in terms of expressions, in 1891 Frege emphasized the distinction between an expression and its reference, and made it clear that functions are what are \emph{denoted} by the corresponding  syntactic expressions \cite{Frege91}. He went on to note that his conception of function extends previous ones in two ways: by enlarging the collections of signs that could be used to construct a functional expression, and by enlarging the domain of possible arguments for functions.\footnote{See \cite[pp.~137 and 140]{FregeReader}.} Regarding the first extension, Frege allowed signs such as the equality symbol to occur in functional expressions, thus allowing ``$x^3 + x^2 + x + 1 = 0$'' to be classified as such. Regarding the second extension, Frege wrote:
\begin{quote}
 Not merely numbers, but objects in general, are now admissible; and here persons must assuredly be counted as objects.\footnote{``Es sind nicht mehr blo{\ss} Zahlen zuzulassen, sondern Gegenst\"{a}nde \"{u}berhaupt, wobei ich allerdings auch Personen zu den Gegenst\"{a}nden rechnen mu\ss.''} \cite[p.~17]{Frege91}
\end{quote}
The syntax of Frege's logical language allows, moreover, for higher-order functionals, which is to say, functions which take elements of a domain of functions as arguments. This feature prompted the following observation in the \emph{Begriffsschrift}:
\begin{quote}
One sees here particularly clearly that the concept of function in Analysis, which in general I have followed, is far more restricted than that developed here.\footnote{``Man sich hieran besonders klar, dass der Functionsbegriff der Analysis, dem ich mich im Allgemeinen angeschlossen habe, weit beschr\"ankter ist als der hier entwickelte.''} \cite[\S 10]{Frege1879} 
\end{quote}
Although Frege was influenced by Riemann and developments in the theory of complex functions \cite{tappenden:06}, there are features of his treatment of the function concept that distinguish it from the contemporary mathematical notion. Nonetheless, in the context of the history we have sketched here, Frege's dramatic expansion was novel and bold.

Today, we are apt to look back at the history and wonder what mathematicians before Dedekind and Frege were missing, that is, why they didn't realize the extent to which they could simplify terminology and notation and eliminate conceptual clutter. After all, all they had to do was write $f : A \to B$ to denote a function $f$ between two arbitrary domains, $A$ and $B$, and recognize sequences, number-theoretic functions, permutations, and transformations as instances of such. But this seems to us to be the wrong question to ask. A better one is this: why would anyone want to view these patently different entities as instances of a single concept? In practice, they are described in very different ways: analytic functions were typically given by piecewise analytic expressions; number theoretic functions were given by implicit algorithms, such as counting the numbers less than and relatively prime to some number; Galois' substitutions were given by explicit lists; early geometric transformations were given by geometric constructions. Moreover, the different kinds of objects support very different operations: functions in analysis can sometimes be composed, but sequences and series cannot; substitutions (and Cayley's ``operations'') can always be inverted, whereas number-theoretic functions generally cannot; one can sum an expression over all the permutations of a finite set, while one certainly cannot sum an expression over all the functions from $\RR$ to $\RR$; notions of continuity and differentiability make no sense outside the realm of analysis; and so on. In a discussion of the function concept in analysis, Renaud Chorlay nicely sums up the state of affairs as follows:
\begin{quote}
As du Bois-Reymond strikingly put it, the ``most general'' function, the ``arbitrary function,'' the ``function on which no hypotheses is made'' is something about which nothing can be said. It is by no means an object to be studied, it is but an (intensionally) empty place in the whole epistemic configuration: not something to investigate, but a kind of background against which ever more specific function classes can be delineated; meaningless (on the intensional level) because all-encompassing (on the extensional level). \cite{chorlay:unp}
\end{quote}
This characterization is even more applicable to the notion of a function between two arbitrary domains.

There are invariably startup costs involved in fashioning coherent means of dealing with new abstract objects, and significant concerns. Thus one should not expect to see an investment in the unification of the function concept until there was either a pressing need or clear benefits to be had. Our study below will suggest that some of these benefits stem from the uniform treatment of algebraic structures composed of functions. For example, the set of automorphisms of a structure form a group under composition; but there is little reason to recognize permutations, geometric transformations, or automorphisms of a field as such until one has particularly useful things to say about groups of automorphisms in general. Similarly, one can view a number of common mathematical constructions as quotienting by the kernel of a suitable function, but there is no point to making the effort to recognize them as such until one realizes that substantial things can be said about the commonality. Alternatively, it is also possible that the unification of the function concept is best viewed as a by-product of the drive towards subsuming mathematical objects under a single unifying foundation, a move which, in turn, might have been encouraged by other mathematical needs.

It is beyond the scope of our study to speculate more than this as to the factors that encouraged the development of a general concept of function and led to its gradual acceptance. But noting the absence of an overarching function concept in the early nineteenth century serves to clarify the nature of the project here. When we consider the history of number-theoretic characters in proofs of Dirichlet's theorem, we are studying the evolving mathematical treatment of \emph{what we now take to be} certain kinds of functions, in the hopes that understanding the forces that guided that evolution will help illuminate the forces that shaped the evolution of modern mathematics more generally. In doing so, we will focus on the way characters were used and the properties that were ascribed to them, while, for the most part, steering clear of the question as to how and why they came to be viewed as instances of the general function concept.

\subsection{The evolution of the function concept}

We have seen that it was not until the latter half of the nineteenth century that there was any hint of a general notion of a function between any two domains, and that it was a while before what we currently take to be instances of the function concept were subsumed under such a general notion. Let us label this trend the \emph{generalization} (or, perhaps, \emph{unification}) of the function concept. The point of this section is to note that this is not the only sense in which the function concept evolved over the course of the century.

To start with, there was also the \emph{extensionalization} of the function concept. Roughly speaking, by an ``extensional'' view we mean the implicit understanding that a function is completely determined by its extension, that is, the values it takes at all the inputs in its domain. In contrast, by an ``intensional'' view of the function concept we mean any way of thinking that presupposes that there is more to a function than its bare extension. For example, let $f$ be the function from real numbers to real numbers given by $f(x) = x^2 + 2 x + 1$. Saying that ``$f$ has three terms'' or that ``the second term of $f$ is $2 x$'' amounts to treating $f$ intensionally, because these statements cannot be couched solely in terms of the values that $f$ takes. Of course, most mathematicians today would acknowledge that such phrases refer to ``the expression of $f$'' or ``the definition of $f$'' or ``the representation of $f$,'' rather than the function $f$ itself. But one of the themes that will emerge below is that this distinction was not clear in the nineteenth century, and fundamental considerations pulled in opposite directions. For concreteness, it is natural to think of functions as expressions, or, at least, entities somehow associated with a concrete representation. But, as we will see, there are also advantages to thinking of operations on functions that depend only on the abstract input-output relation. Our case study will highlight an uneasy tension between these two ways of thinking.

Recall that in Euler's time, a function from the real numbers to the real numbers was expected to have a representation in terms of a certain kind of analytic expression. Arguments involving functions could then make use of such a representation. Precisely this feature of Cauchy's treatment of Fourier series was the target of Dirichlet's criticism in the opening paragraphs of his 1829 paper \cite{Dirichlet1829} on that same subject. Cauchy, like Dirichlet, was concerned with the question as to when the Fourier series of a real-valued function converges to the value of the original function at a given point. But Cauchy's analysis presupposed that the function in question is represented by a power series, which could then be applied to complex arguments. The use of the phrase ``arbitrary function'' in Dirichlet's title signaled his intention to lift this restriction. Even though the paper dealt primarily with continuous functions, the hypothesis of continuity was expressed purely extensionally, which is to say, in terms of the values that a function takes at its arguments. In other words, Dirichlet took great pains not to treat analytic operations like differentiation as symbolic operations, but, rather, as operations on the functions themselves, viewed extensionally. As Dirichlet famously emphasized in that paper, this renders the concept of function open-ended, in that all that is needed is some determinate relationship between input and output values. In particular, this made it possible for him to consider, as a difficult example, functions which take one value on the rationals, and another value on the irrational numbers.

This brings us to another aspect of the evolution of the notion of a function, which we will call the \emph{liberalization} of the function concept. Setting aside issues of unification and generality, there is also the issue of the language and methods one helps oneself to in \emph{defining} functions in a particular domain. As Frege later noted,\footnote{See the discussion in Section~\ref{analysis:section}.} the very act of accepting a definition by cases depending on whether an argument is rational or not was a striking move on Dirichlet's part. This paves the way to the use of more dramatically non-constructive methods in describing functions from the reals to the reals, say, in terms of limits and other infinitary operations. Riemann's contributions to function theory were especially novel in this regard. Some, like Dedekind, hailed the fact that Riemann's methods made it possible to consider functions without having a particular representation to work with, or method of computation. Others, like Weierstrass, held the Riemannian approach to be defective for just this reason.\footnote{See Bottazzini and Gray \cite{bottazzini:gray:13}. For an interesting exploration of the ways that nineteenth century analysis expanded to incorporate a more liberal understanding of the function concept, see Chorlay \cite{chorlay:unp}.}

Finally, there is the issue of the \emph{reification} of the function concept, which is to say, the decision to treat functions as objects, on par with mathematical objects like the natural numbers. In Section~\ref{objects:section}, we will expand on what this amounts to.

We have already indicated that we will have little more to say about the generalization of the function concept.
And, since characters are fairly simple combinatorial objects, we can avoid issues having to do with liberal uses of the infinite. We will be quite concerned, however, with issues regarding the reification of characters, and their treatment as extensional objects. Thus our case study considers some important aspects of the evolution of the function concept, but ignores others, as well as the broader question as to how all the various components are related to one another.

\section{Dirichlet's theorem}
\label{contemporary:section}

Two integers, $m$ and $k$, are said to be \emph{relatively prime}, or \emph{coprime}, if they have no common factor. In 1837, Dirichlet proved the following:
\begin{theorem}
\label{dirichlet:theorem}
 If $m$ and $k$ are relatively prime, the arithmetic progression $m, m + k, m + 2k, \ldots$ contains infinitely many primes.
\end{theorem}
In other words, if $m$ and $k$ are relatively prime, there are infinitely many primes congruent to $m$ modulo $k$. Dirichlet pointed out that Legendre assumed this fact, without proof, in 1788 \cite{Legendre88}, when proving the law of quadratic reciprocity. In his \emph{Disquistiones Arithmeticae}, Gauss noted this gap in Legendre's work after presenting his own proof of the law of quadratic reciprocity, one which does not rely on Theorem~\ref{dirichlet:theorem}.\footnote{See Gauss' remarks in \cite[\S 150]{gauss:01}.} But Gauss himself was never able to prove that theorem. Dirichlet's own proof is striking, not only due to the fact that it finally established Legendre's conjecture, but also due to its sophisticated  use of the methods of analysis in establishing a purely number-theoretic assertion. Dirichlet noted \cite[pp.~309--310]{Dirichletshort} that his method was inspired by a proof that there are infinitely many primes due to Euler \cite[Chapter XV]{Euler48}, though Dirichlet's ideas go considerably beyond Euler's. In this section, we will describe Euler's proof, and then define the modern notion of a group-theoretic character, which supports the generalization to arbitrary characters. We will then describe contemporary proofs of Dirichlet's theorem, using modern terminology and notation, as presented in textbooks such as Everest and Ward \cite[pp.~207--224]{EverWard05}. In Section~\ref{objects:section}, we will use this presentation to frame our discussion of the historical development, which occurs in Sections~\ref{dirichlet:section} and \ref{evolution:section}.

There is a sense, however, in which our modern presentation is old-fashioned. The use of $m$ and $k$ for the initial value and common difference in the statement of Theorem~\ref{dirichlet:theorem} is due to Dirichlet, and was picked up by a number of successive authors, including Dedekind and Hadamard. Modern presentations are more apt to use $a$ and $d$, but although the use of $m$ and $k$ may feel alien to readers familiar with contemporary textbook proofs, it will facilitate the historical comparisons later on. (Our reason for using the variable $q$ to range over the prime numbers will similarly become clear when we discuss Dirichlet's original proof.)

\subsection{Euler's proof that there are infinitely many primes}

In the \emph{Elements}, Euclid proved that there are infinitely many primes, but his proof does not provide much information about how they are distributed. Euler, in his \emph{Introductio in Analysin Infinitorum} \cite{Euler48}, proved the following:
\begin{theorem}
\label{euler:thm}
The series $\sum_{q}\frac{1}{q}$ diverges, where the sum is over all primes $q$.
\end{theorem}
This implies that there are infinitely many primes, but also says something more about their density. For example, since we know that the series $\sum_{n}\frac{1}{n^{2}}$ is convergent, it tells us that, in a sense, there are ``more'' primes than there are squares.

Euler's proof of Theorem~\ref{euler:thm} centers around his famous zeta function,
\[
\zeta(s) = \sum_{n=1}^{\infty}n^{-s},
\]
defined for a real variable $s$. (The zeta function was later extended by Riemann to the entire complex plane via analytic continuation.) It is not hard to show that the series $\zeta(s)$ converges uniformly on the interval $[a, \infty)$, where $a$ is any number strictly greater than 1. For $s > 1$, the infinite sum can also be expressed as an infinite product:
\[
 \sum_{n=1}^{\infty} n^{-s}= \prod_{q}\left(1-\frac{1}{q^{s}}\right)^{-1},
\]
where the product is over all primes $q$. This is known as the \emph{Euler product formula.} Roughly, this holds because we can write each term of the product as the sum of a geometric series,
\[
\left(1-\frac{1}{q^{s}}\right)^{-1} = 1 + q^{-s} + q^{-2s} + \ldots
\]
and then expand the product into a sum. The unique factorization theorem tells us that every integer $n > 1$ can be written as a product $q_1^{i_1} \cdot q_2^{i_2} \cdots q_k^{i_k}$, and so the term $n^{-s} = q_1^{-i_1 s} \cdot q_2^{-i_2 s} \cdots q_k^{-i_k s}$ will occur exactly once in the expansion. Since we are dealing with infinite sums and products, the Euler product formula implicitly makes a statement about limits, and some care is necessary to make the argument precise; but this is not hard to do.

If we take the logarithm of each side of the product formula and appeal to properties of the logarithm function, we obtain
\[
 \log\sum_{n=1}^{\infty}n^{-s}=\sum_{q}-\log\left(1-\frac{1}{q^{s}}\right).
\]
Using the Taylor series expansion
\[
\log(1 - x) = -x - \frac{x^2}{2} - \frac{x^3}{3} - \ldots
\]
and changing the order of summations yields
\[
\log\sum_{n=1}^{\infty}n^{-s}=\sum_{q}\frac{1}{q^{s}} + \sum_{n=2}^{\infty}\frac{1}{n}\sum_{q}\frac{1}{q^{ns}}.
\]
At this stage, keep in mind that we want to show that $\sum_{q}\frac{1}{q}$ diverges, and notice that the first term on the right hand side of the above equation is $\sum_{q}\frac{1}{q^{s}}$. Thus we should consider what happens as $s$ tends to 1 from above. It is not hard to show that the second term on the right-hand side is bounded by a constant that is independent of $s$, a fact that can be expressed using ``big O'' notation as follows:
\begin{equation}
\label{euler:primes:eqn}
\log\sum_{n=1}^{\infty}n^{-s}=\sum_{q}\frac{1}{q^{s}} + O(1).
\end{equation}
As $s$ approaches $1$ from above, the left-hand side clearly tends to infinity. Thus, the right-hand side, $\sum_{q}\frac{1}{q^{s}}$, must also tend to infinity, which implies that $\sum_{q}\frac{1}{q}$ diverges.

\subsection{Group characters}
\label{group:character:section}

Just as Euler's proof shows that there are infinitely many primes by establishing that the series $\sum_{q}\frac{1}{q}$ is divergent, proofs of Dirichlet's theorem establish that there are infinitely many primes $q$ congruent to $m$ modulo $k$ by showing that the series $\sum_{q \equiv m \pmod{k}}\frac{1}{q}$ is divergent. Remember that we are assuming that $m$ and $k$ are relatively prime. In general, the residues modulo $k$ that are relatively prime to $k$ form a multiplicative group. To adapt the argument above, we need series that are more refined than the zeta function, and a device that enables us to focus attention on the integers that are congruent to $m$ modulo $k$. This is where the notion of a group-theoretic character comes in. 

Notice that the expression $\sum_{q \equiv m \pmod{k}}\frac{1}{q}$ can be written $\sum_q \frac{f(q)}{q}$, where $f(n)$ is equal to $1$ if $n$ is congruent to $m$ modulo $k$ and $0$ otherwise. With this observation, it is natural to try to emulate Euler's argument with $\sum_{n = 1} f(n)n^{-s}$ in place of the $\zeta$ function, $\sum_{n = 1} n^{-s}$. We will see below that an analogue of the Euler product formula holds if $f$ is \emph{completely multiplicative}, which is to say, it satisfies $f(uv) = f(u)f(v)$ for every $u$ and $v$. The problem is that the particular $f$ just described does not have this property. But Corollary~\ref{orthocorrol} below shows that we can decompose $f$ into a sum of functions that do have this property, and apply the Euler product formula to each component.

Let $G$ be a finite abelian group, written multiplicatively, with an identity denoted by $1$. A \emph{character $\chi$ on $G$} is a homomorphism from $G$ to the nonzero complex numbers, $\CC^*$. In other words, a character $\chi$ is a nonzero function satisfying $\chi(g_1 g_2) = \chi(g_1) \chi(g_2)$ for every $g_1$ and $g_2$ in $G$. Since $\chi(1) = \chi(1 \cdot 1) = \chi(1) \chi(1)$ and $\chi(1)$ is nonzero, we have $\chi(1) = 1$. Since $G$ is a finite group, for every $g$ there is a least $n > 1$ satisfying $g^n = 1$; this $n$ is called the \emph{order} of $g$ and denoted $o(g)$. The fact that $\chi(g)^{o(g)} = \chi(g^{o(g)}) = \chi(1) = 1$ means that $\chi(g)$ is a ``root of unity'' for every $g$ in $G$. The character which is equal to $1$ for every $g$ in $G$ is called the \emph{trivial character} and denoted by $\chi_0$.

Define the product $\chi \cdot \psi$ of two characters pointwise, by
\[
 (\chi \cdot \psi)(g) = \chi(g) \psi(g),
\] 
for every $g$ in $G$. This multiplication is commutative, and we have $\chi \cdot \chi_0 = \chi$ for every character $\chi$, which is to say, $\chi_0$ is a multiplicative identity. Recall that if $\omega$ is any complex root of unity, its complex conjugate, $\overline \omega$, is also a root of unity, satisfying $\omega \overline \omega = 1$. We can also lift the operation of conjugation to characters, defining $\overline \chi$ by the equation $\overline \chi(g) = \overline{\chi(g)}$ for each $g$. Then clearly we have $\chi \cdot \overline \chi = \chi_0$. In other words, the set of characters on $G$ forms an abelian group, with $\chi^{-1} = \overline \chi$. We will denote this group $\widehat{G}$.

The following theorem is fundamental:
\begin{theorem}
\label{GisoGhat}
If $G$ is any finite abelian group, $\widehat{G}$ is isomorphic to $G$. In particular, $|\widehat{G}| = |G|$.
\end{theorem}

To see this, first consider the case where $G = \langle g \rangle = \{1, g, g^2, \ldots, g^{n-1}\}$ is a cyclic group generated by an element $g$ of order $n$. Then any character $\chi$ on $G$ has to map $g$ to an $n$th root of unity, $\omega$. This determines the behavior of $\chi$ completely, since then we have $\chi(g^i) = \omega^i$ for every $i$. Thus we can let $\chi_\omega$ denote the unique character that maps $g$ to $\omega$. But now if we let $\omega = e^{2 \pi i / n}$, then $\omega$ is what is known as a ``primitive root of unity,'' which is to say that all the roots of unity are given by $1, \omega, \omega^2, \ldots, \omega^{n-1}$. Notice that these roots form a multiplicative group that is isomorphic to $G$. It is easy to verify that the map which sends the element $g^i$ of $G$ to the character $\chi_{\omega^i}$ of $\widehat G$ is an isomorphism.

In the more general case, we appeal to the structure theorem for finite abelian groups, which says that any such group $G$ can be written as a product $G_1 \times \ldots \times G_k$ of cyclic groups. This means that every element $g$ of $G$ can be written uniquely as a product $g = g_1 g_2 \cdots g_k$, where each $g_i$ is in $G_i$. Given characters $\chi_1, \chi_2, \ldots, \chi_k$ on $G_1, G_2, \ldots, G_k$ respectively, one gets a character $\chi$ on $G$ defined by $\chi(g) = \chi_1(g_1) \chi_2(g_2) \cdots \chi_k(g_k)$, where $g = g_1 g_2 \ldots g_k$ is the decomposition described above. Moreover, it is not hard to show that every character on $G$ arises in this way. This shows that $\hat G$ is isomorphic to $\hat G_1 \times \hat G_2 \times \ldots \times \hat G_k$. By the analysis of the cyclic case, the latter is, in turn, isomorphic to $G_1 \times G_2 \times \ldots \times G_k$, and hence to $G$.

The following theorem expresses two important properties, known as the ``orthogonality relations'' for group characters.

\begin{theorem}
\label{ortho}
Let $G$ be finite abelian group.  Then for any character $\chi$ in $\widehat{G}$, we have
\[
\sum_{g \in G}\chi(g) =
\begin{cases}
  |G| & \mbox{if $\chi = \chi_{0}$} \\
  0   & \mbox{if $\chi \neq \chi_{0}$,}
\end{cases}
\]
and for any element $g$ of $G$, we have
\[
\sum_{\chi \in \widehat{G}}\chi(g) =
\begin{cases}
  |G| & \mbox{if $g=1_{G}$} \\
  0   & \mbox{if $g \neq 1_{G}$.}
\end{cases}
\]
\end{theorem}

The first equation clearly holds when $\chi$ is the trivial character, $\chi_0$, since, in this case, each term of the sum is equal to $1$. Otherwise, pick $h$ such that $\chi(h) \neq 1$ and note
\[
 \chi(h) \sum_{g \in G} \chi(g) = \sum_{g \in G} \chi(h g) = \sum_{g \in G} \chi(g),
\]
since $hg$ ranges over the elements of $G$ as $g$ does. Since $\chi(h) \neq 1$, we must have $\sum_{g \in G} \chi(g) = 0$. The second equation can be established in a similar way.

The orthogonality relations make it possible to do ``finite Fourier analysis,'' in the following sense: if $f$ is any function from $G$ to the complex numbers and we define the ``Fourier transform'' $\hat f$ of $f$ by $\hat f(\chi) = \sum_g f(g) \chi(g)$, then $f$ can be recovered from its Fourier transform: $f = \frac{1}{|G|} \sum_\chi \hat f (\chi) \chi$. This shows, in particular, that any function from $G$ to the complex numbers can be written as a linear combination of characters. The second orthogonality relation also provides the following useful corollary:

\begin{corollary}
\label{orthocorrol}
For any $g, h \ \in \ G$ we have the following:
\[
\sum_{\chi \ \in \ \widehat{G}}\chi(g)\overline{\chi(h)}= \begin{cases} |G| & \ \mbox{if} \ g=h \\ 0 & \ \mbox{if} \ g \neq h\end{cases}
\]
\end{corollary}

This follows from the fact that we have
\[
\sum_{\chi \ \in \ \widehat{G}}\chi(g)\overline{\chi(h)} =\sum_{\chi \ \in \ \widehat{G}}\chi(g)\chi(h)^{-1}=\sum_{\chi \ \in \ \widehat{G}}\chi(gh^{-1})= \begin{cases} |G| & \ \mbox{if } \ g=h \\ 0 & \ \mbox{if } \ g \neq h \end{cases}.
\]
This corollary will enable us to focus on the residue class of $m$ modulo $k$ in the proof of Dirichlet's theorem.

\subsection{Dirichlet characters and $L$-series}
\label{dirichlet:characters:section}

Let $k$ be an integer greater than or equal to $1$. It is a fundamental theorem of number theory that an integer $n$ is relatively prime to $k$ if and only if $n$ has a multiplicative inverse modulo $k$; in other words, if and only if there is some $n'$ such that $n n' \equiv 1 \bmod{k}$. This implies that the residue classes of integers modulo $k$ that are relatively prime to $k$ form a group, denoted $(\ZZ/k\ZZ)^*$, with multiplication modulo $k$. The cardinality of $(\ZZ/k\ZZ)^*$, that is, the number of residues relatively prime to $k$, is denoted $\ph(k)$, and the function $\ph$ is called the \emph{Euler phi function}.

Any character $\chi$ on $(\ZZ/k\ZZ)^*$ can be ``lifted'' to a function $X$ from $\ZZ$ to $\CC$ defined by
\[
 X(n) =
  \begin{cases}
    \chi(n \bmod k) &\text{if $n$ is relatively prime to $k$} \\
    0 &\text{otherwise.}
  \end{cases}
\]
Such a function is called a \emph{Dirichlet character modulo $k$}. Dirichlet characters are completely multiplicative, which is to say, $X(mn) = X(m)X(n)$ for every $m$ and $n$ in $\ZZ$. Mathematicians typically use the symbol $\chi$ to range over Dirichlet characters, blurring the distinction between such functions and their group-character counterparts. This is harmless, since there is a one-to-one correspondence between the two, and so we will adopt this practice as well.

We can now generalize the method of Euler's proof. Roughly speaking, we need a variant of the zeta function that will allow us to focus on primes in a particular residue class modulo $k$. To that end, given a Dirichlet character $\chi$ modulo $k$, define the \emph{Dirichlet $L$-function}, or \emph{$L$-series},
\[
  L(s, \chi)= \sum_{n=1}^{\infty}\frac{\chi(n)}{n^{s}},
\]
where $s$ is any complex number. This formal series will converge whenever $\mathfrak{R}(s) > 1$, that is, the real part of $s$ is greater than one.
And just as $\sum_{n=1}^{\infty}\frac{1}{q^{s}}$ can be written as a product via the Euler product formula, so each $L(s, \chi)= \sum_{n=1}^{\infty}\frac{\chi(n)}{n^{s}}$ has a useful product expansion.
\begin{theorem}
\label{modernEP}
Let $\chi$ be a Dirichlet character modulo $k$.  Then the $L$-function associated with $\chi$ has an Euler product expansion for $\mathfrak{R}(s)>1$,
\[
L(s, \chi) = \sum_{n=1}^{\infty}\frac{\chi(n)}{n^{s}}=\prod_{q}\left(1-\frac{\chi(q)}{q^{s}}\right)^{-1}=\prod_{q\nmid k} \left(1-\frac{\chi(q)}{q^{s}}\right)^{-1}.
\]
\end{theorem}
The last identity follows from the fact that, in the product, we can ignore those primes $q$ that divide $k$, since $\chi(q) = 0$ for such $q$. With the product formula in place, we can sketch a proof of Dirichlet's theorem.

\subsection{Contemporary proofs of Dirichlet's theorem}
\label{contemporary:proofs:section}

Recall that we want to prove that there are infinitely many primes $q$ such that $q\equiv m \pmod{k}$, where $m$ and $k$ are relatively prime.  As in the proof that $\sum_{q}\frac{1}{q^{s}}$ diverges, we begin by taking logarithms of both sides of the Euler product expansion for $L(s, \chi)$, where $\chi$ is a Dirichlet character modulo $q$:
\[
\log L(s, \chi)= - \sum_{q \nmid k}\log\left(1 - \frac{\chi(q)}{q^{s}}\right).
\]
As before, we make use of the Taylor series expansion for the logarithm on the right hand side of the above equation to obtain:
\begin{align*}
\log L(s, \chi)
& = \sum_{q \nmid k}\sum_{j=1}^{\infty}\frac{1}{j}\frac{\chi(q^{j})}{q^{sj}} \\
& = \sum_{q \nmid k} \frac{\chi(q)}{q^{s}} + \sum_{q \nmid k, j=2}^\infty \frac{1}{j}\frac{\chi(q^{j})}{q^{sj}}.
\end{align*}
One can show that the second term in the expression is bounded by a constant that is independent of $s$ and $\chi$, which can be expressed as follows:
\[
\log L(s, \chi) = \sum_{q \nmid k}\frac{\chi(q)}{q^{s}} \ + \ O(1).
\]
Now comes the crucial use of Corollary~\ref{orthocorrol} to pick out the primes in the relevant residue class.  We multiply each side of the above equation by $\overline{\chi(m)}$ and then take the sum of these over all the Dirichlet characters modulo $k$. (Recall that we can identify each Dirichlet character with the corresponding group character, that is, the corresponding element of $\widehat{(\mathbb{Z}/k \mathbb{Z})^*}$.) Thus we have:
\begin{equation*}
\sum_{\chi \in \widehat{(\mathbb{Z}/k \mathbb{Z})^*}} \overline{\chi(m)}\log L(s, \chi) =
\sum_{\chi \in \widehat{(\mathbb{Z}/k \mathbb{Z})^*}} \overline{\chi(m)}\sum_{q \nmid k}\frac{\chi(q)}{q^{s}}\  + \  O(1).
\end{equation*}
To simplify this expression, we exchange the summations on the right-hand side, and appeal to Corollary~\ref{orthocorrol}. Since the cardinality of the group $(\mathbb{Z} / k \mathbb{Z})^*$ is $\ph(k)$, we obtain
\begin{equation}
\label{sumoverchar2}
\sum_{\chi \in \widehat{(\mathbb{Z} / k \mathbb{Z})^*}}\overline{\chi(m)}\log L(s, \chi) =
\ph(k) \sum_{q \equiv m \pmod{k}}\frac{1}{q^{s}} \ + \ O(1).
\end{equation}
This is analogous to the equation (\ref{euler:primes:eqn}) in Euler's proof. Our goal is once again to show that the left-hand side tends to infinity as $s$ approaches 1 from above; this implies that the right-hand side tend to infinity, which, in turn, implies that there infinitely many primes $q$ that are congruent to $m$ modulo $k$. However, now the left-hand side is considerably more complicated than the expression $\log \sum_{n=1}^{\infty} n^{-s}$ in Euler's proof.

To show that $\sum_{\chi \in \widehat{(\mathbb{Z} / k \mathbb{Z})^*}}\overline{\chi(m)}\log L(s, \chi)$ tends to infinity a $s$ approaches $1$, we divide the characters into three classes, as follows:
\begin{enumerate}
\item The first class contains only the principal character $\chi_0$, which takes the value of 1 for all arguments that are relatively prime to $k$, and 0 otherwise.
\item  The second class consists of all those characters which take only real values (i.e.\ 0 or $\pm 1$), other than the principal character.
\item  The third class consists of those characters which take at least one complex value.
\end{enumerate}
It is not difficult to show that $L(s, \chi_{0})$ has a simple pole at $s=1$, which implies that the term $\overline{\chi_0}(m) \log L(s,\chi_0)$ approaches infinity as $s$ approaches $1$. The real work involves showing that for all the other characters $\chi$, $L(s,\chi)$ has a finite nonzero limit. This implies that the other terms in the sum approach a finite limit, and so the entire sum approaches infinity.

For characters in the third class, that is, the characters that take on at least one complex value, the result is not difficult. For characters in the second class, the result is much harder, and Dirichlet used deep techniques from the theory of quadratic forms to obtain it. In the years that followed, other mathematicians found alternative, and simpler, ways of handling this case. But even in modern presentations, this case remains the most substantial and technically involved part of the proof.

Our presentation has been thoroughly ``modern.'' In the next section, we will consider some of the methodological features of the proof that make it so. This will enable us to draw interesting contrasts with Dirichlet's proof, and then explore the way that presentations of Dirichlet's theorem gradually took on such a modern character.

\section{Modern aspects of contemporary proofs}
\label{objects:section}


Although our presentation uses contemporary terminology and notation, there is a sense in which it is a faithful description of Dirichlet's 1837 proof. Dirichlet did not, and could not, rely on a general notion of group character, as the general notion of a group was not articulated before Cayley did so in 1854 \cite{cayley:54}, and was not brought into general currency until Kronecker's 1870 paper \cite{kronecker:70}, which first presented the structure theorem for finite abelian groups.\footnote{See Wussing \cite{wussing:84} for the history of group theory.} Another important difference is that even though Dirichlet's argument used complex numbers in a central way, there was less established background in complex analysis than is available today, and Dirichlet tended to reduce the calculations to real analysis whenever possible. Thus his variable $s$ ranged over real numbers, and his calculations involve real-valued sines and cosines where today we are comfortable sticking with the complex exponential. In addition, we have already noted that many of the technical details were streamlined over the years. Despite all this, the outline above characterizes the central ideas of his proof, and most mathematicians would not find it unreasonable to say that that is, essentially, how he obtained the result.

But, as we will see in Section~\ref{dirichlet:section}, there is one very striking difference: in Dirichlet's original presentations, there is no mention of characters at all. That is, Dirichlet's papers contain certain expressions that we now recognize as values of the various characters, and summations that are \emph{tantamount} to summing over all the characters. But the characters themselves are only objects that we project back into the argument from our current understanding. They hover over the page as shade-like premonitions, ghosts of mathematics yet to come.

In Section~\ref{evolution:section}, we will consider presentations of Dirichlet's theorem given by Dedekind, de la Vall\'ee-Poussin, Hadamard, Kronecker, and Landau, and see how the characters were gradually brought to life. We will see that many of the benefits of giving characters a substantive embodiment are notational and pragmatic, but that is not to say that they are \emph{merely} notational and pragmatic: treating characters as bona fide objects comes with serious mathematical constraints and obligations, and provides conceptual reorientations that have great bearing on the kind of mathematics we do, and the way we do it. We will argue that the reification of the notion of a character is a prototypical instance of the conceptual changes that are hallmarks of the transition to ``modern'' mathematical thought, and that understanding how and why the changes came about shed light on the way we do mathematics today.

But the observations above present us with a terminological conundrum: should we describe the various historical texts as ``versions'' or ``presentations'' of Dirichlet's proof, or different proofs entirely? Having raised this issue, we will, for the most part, set it aside, and be fairly cavalier with our terminology. Since our specific concern is to study the way that language, conceptualization, and inferential practice evolved over the years, and the effects that had on the mathematics, we need not explicitly address the question as to when it is proper to consider two proofs essentially the same or essentially different.

Talking about historical texts in modern terms is difficult, and it is always misleading to portray the history of mathematics as a muddled and inefficient attempt to arrive at the contemporary enlightened view. We hope we have not fallen into this trap. If there is anything that deserves to be treated as a rational pursuit, mathematics should count as such, and so one would expect there to be good reasons that we do mathematics the way we do. At the same time, there are often good reasons to question the way we do mathematics today, and recognize there are tradeoffs involved in the historical decisions that were made. Comparing mathematical texts from different historical eras and trying to understand what has changed provides a fruitful way of understanding the values that drive mathematical change. But it is often easier to draw contrasts by starting with the mathematics with which we are most familiar, and so, having discussed contemporary approaches to Dirichlet's theorem, let us foreshadow some of the contrasts we wish to consider.


In Section~\ref{functions:section}, we noted that characters (whether we refer to group characters, or Dirichlet characters) are instances of the contemporary function concept, a concept which evolved significantly over the course of the nineteenth century. These are some of the salient features of the treatment of characters in our modern presentation:
\begin{enumerate}
\item Group characters are given an abstract, axiomatic definition as functions that satisfy the homomorphism property, and the Dirichlet characters are introduced as a natural extension of the notion.
\item In particular, one defines the set of characters modulo $k$ extensionally. Only later does one show that this set is finite, and provide explicit ways of describing and enumerating them.
\item Characters are studied in their own right, and their general properties are enunciated in propositions and theorems.
\item One sums over sets of characters, without needing representations for any particular one. More, generally, one characterizes operations on characters (such as the product of two characters) extensionally, and not in terms of their representations.
\item Characters appear as arguments to other functions, namely, the $L$-functions. In particular, it is clear from the definition that $L(s,\chi)$ depends only on the extension of $\chi$.
\item One defines various sets of characters extensionally, for example, distinguishing the trivial, real, and complex characters in terms of the values they assume. More generally, one typically carries out arguments without making reference to any particular representation.
\item The characters modulo $k$ are viewed as elements of an algebraic structure, namely, a group, with multiplication defined pointwise and the trivial character serving as identity.
\end{enumerate}
In Section~\ref{dirichlet:section}, we will discuss Dirichlet's original proof, and in Section~\ref{evolution:section}, we will consider the way the proof was gradually transformed to reflect our contemporary understanding. We will see that proofs along the way possess various subsets of the properties just enumerated, and that, in some cases, the authors are noticeably squeamish, or at least self-conscious, of these features. Importantly, Dirichlet's original proof has \emph{none} of the properties just enumerated, providing a clear contrast to the style of presentation that is common today.

One way of characterizing the difference between Dirichlet's original proof and our modern presentation is to say that, in the latter, characters are treated as full-fledged mathematical objects, whereas there are no such objects in Dirichlet's version. Elsewhere \cite{avigad:morris:b} we reflect in greater detail on what it means to say that a piece of mathematics sanctions certain entities as ``objects.'' Here, let us merely summarize some of the senses in which this can be said to be the case in the modern proofs of Dirichlet's theorem:
\begin{enumerate}
 \item Characters fall under a recognized grammatical category, which allows us to state things about them and define operations and predicates on them.
 \item There is a clear understanding of what it means for two expressions to represent the \emph{same} character, namely, that they take the same values at all arguments. Linguistic conventions ensure that expressions occurring in a proof respect this ``sameness.'' For example, the expression $L(s, \chi)$ does not depend on the data used to represent $\chi$. 
 \item One can quantify and sum over characters; in logical terms, they can fall under the range of a bound variable.
 \item One can define functions which take characters as arguments, and sets of characters; indeed, characters can be elements of arbitrary algebraic structures.
\end{enumerate}
What these features have in common is that they are fundamental to the way we do mathematics, bearing upon the proper use of mathematical language and inference at a very low level.

\section{Dirichlet's original proof}
\label{dirichlet:section}

We have already asserted that characters appear only ``implicitly'' in Dirichlet's original proof \cite{Dirichlet37}. There is nothing mysterious about this: what we mean is that, in Dirichlet's proof, there are certain symbolic expressions that we now recognize as denoting the values of characters; and that, moreover, some of Dirichlet's calculations and inferences invoke what we now recognize as general properties of characters.

Let us spell out the details. Like Dirichlet, we will first consider the case were the common difference in the arithmetic progression is a prime number, denoted by $p$ instead of $k$. It was well known in the nineteenth century that one can always find a \emph{primitive element} modulo $p$, which is to say, an element $c$, such that the $p-1$ residues $c^0, c^1, c^2, \ldots, c^{p-2}$ modulo $p$ yields all the nonzero resides, $1, 2, \ldots, p - 1$ (not necessarily in the same order). In modern terms, we would say that the group $(\ZZ / p \ZZ)^*$ of units modulo $p$ is cyclic, generated by the residue class of $c$. In more elementary terms, this amounts to saying that for every integer $n$, there is a number $\gamma_n$ with the property that $c^{\gamma_n} \equiv n \bmod p$. We saw in Section~\ref{group:character:section} that if $\chi$ is a Dirichlet character modulo $p$ (that is, $\chi$ corresponds to  character on $(\ZZ / p \ZZ)^*$), then $\chi(c)$ is a $p$th root of unity, say, $\omega$; and, moreover, $\chi$ is entirely determined by $\omega$, in the sense that for every $n$ relatively prime to $p$, $\chi(n) = \omega^{\gamma_n}$. Thus Dirichlet simply wrote $\omega^{\gamma_n}$ where we would write $\chi(n)$. The notation presupposes that one has fixed a choice of the primitive element, $c$, though any primitive element will work equally well.

So far, so good. In the more general case where the modulus is a composite number $k$, however, things get more complicated. First, write $k$ as a product of primes,
\[
 k = 2^\lambda p_1^{\pi_1} p_2^{\pi_2} \cdots p_j^{\pi_j}
\]
where each $p_i$ is an odd prime and $\pi_i$ is greater than or equal to 1. Then the group of units modulo $k$ is isomorphic to the product of the groups of units modulo each term in the factor. Gauss had already shown if $p$ is an odd prime and $\pi$ is an integer greater than or equal to $1$, then one can more generally find a primitive element $c$ modulo $p^\pi$. This means that the residue class of $c$ generates the cyclic group $(\ZZ / p^\pi \ZZ)^*$, or, equivalently, for every $n$ relatively prime to $p$ there is a $\gamma_n$ such that $c^{\gamma_n} \equiv n \bmod p^\pi$. Thus we can choose primitive elements $c_1, \ldots, c_j$ corresponding to $p_1^{\pi_1}, p_2^{\pi_2}, \ldots, p_j^{\pi_j}$. If $\lambda \geq 3$, however, there is no primitive element modulo $2^\lambda$. Rather, $(\ZZ/ 2^\lambda)^*$ is a product of two cyclic groups, and for every $n$ relatively prime to $2^\lambda$ there are an $\alpha_n$ and $\beta_n$ such that $(-1)^{\alpha_n} 5^{\beta_n} \equiv n \bmod 2^\lambda$. Thus for any $n$ relatively prime to $k$, we can write
\[
 n \equiv (-1)^{\alpha_n} 5^{\beta_n} c_1^{\gamma_{1,m}} c_2^{\gamma_{2,m}} \ldots c_j^{\gamma_{j,m}} \bmod k
\]
where each $\gamma_{i,n}$ is the index $n$ relative to $p_i^{\pi_i}$. As above, if we choose appropriate roots of unity $\theta, \ph, \omega_1, \omega_2, \ldots, \omega_j$, we obtain a character
\begin{equation}
\label{characterconstruction}
 \chi(n) = \theta^{\alpha_n} \ph^{\beta_n} \omega_1^{\gamma_{1,n}} \omega_2^{\gamma_{2,n}} \cdots \omega_j^{\gamma_{j,n}}.
\end{equation}
And, once again, every character is obtained in this way. We should note that Dirichlet used the notation $p, p', \ldots$ rather than $p_1, \ldots, p_j$ to denote the sequence of odd primes. Moreover, he used the notation $\alpha, \beta, \gamma, \gamma', \ldots$ to denote the indices, suppressing the dependence on $n$. Thus, Dirichlet wrote $\theta^\alpha \ph^\beta \omega^\gamma \omega'^{\gamma'} \ldots$ for the expression we have denoted $\chi(n)$ above, leaving it up to us to keep in mind that $\alpha, \beta, \ldots$ depend on $n$.

To summarize, in the simple case of a prime modulus $p$, Dirichlet fixed a primitive element modulo $c$, and represented each character $\chi$ in terms of a $p$th root of unity, $\omega$. In that case, the value $\chi(n)$ is given by $\omega^{\gamma_n}$. In the more general case of a composite modulus $k$, Dirichlet fixed primitive elements modulo the terms of the prime factorization of $k$, and represented each character $\chi$ in terms of a sequence $\theta, \ph, \pi, \pi'$ of roots of unity. In that case, the value $\chi(n)$ was written $\theta^\alpha \ph^\beta \omega^\gamma \omega'^{\gamma'} \ldots$, suppressing the information that the exponents $\alpha, \beta, \gamma, \gamma', \ldots$ depend on $n$.

Recall that our contemporary presentation had little to say about \emph{particular} characters, other than the trivial character, $\chi_0$. Rather, characters appear as arguments to the $L$-functions, $L(s,\chi)$, and the proof has us consider summations over the set of all characters. Let us now consider how Dirichlet handled these as well.

Again, with Dirichlet, we begin with the easier case where the common difference of the arithmetic progression is a prime, $p$. Recall that, in that case, each character $\chi$ corresponds to a $p$th root of unity, $\omega$. Dirichlet stated the Euler product formula as follows:
\begin{quote}
We therefore have the equation
\begin{align}
 \prod \frac{1}{1-\omega^{\gamma}\frac{1}{q^{s}}}=\sum\omega^{\gamma}\frac{1}{n^{s}} = L, \label{dirichletEulersimple}
\end{align}
where the multiplication sign ranges over the whole series of primes with the sole exception of $p$, while the summation involves all the integers from $1$ to $\infty$ that are not divisible by $p$. The letter $\gamma$ denotes $\gamma_p$ on the left, and $\gamma_n$ on the right.\footnote{``Man hat daher die Gleichung: \begin{align}
 \prod \frac{1}{1-\omega^{\gamma}\frac{1}{q^{s}}}=\sum\omega^{\gamma}\frac{1}{n^{s}} = L,  \tag{\ref{dirichletEulersimple}}
\end{align} wo sich die Multiplicationszeichen auf die ganze Reihe der Primzahlen, mit alleiniger Ausnahme von $p$, erstreckt, w\"{a}hrend die Summation sich auf alle ganzen Zahlen von 1 bis $\infty$ bezieht, welche nicht durch $p$ teilbar sind.  Der Buchstabe $\gamma$ bedeutet auf der ersten Seite $\gamma_{q}$, auf der zweiten dagegen $\gamma_{n}$.'' \cite[pp.~317--318]{Dirichlet37} We have replaced Dirichlet's equation number with our own, and throughout this section we have modified the translation cited in \cite{Dirichlet37}.} \cite[p.~3]{Dirichlet37}
\end{quote}
Compare this to the statement of Theorem~\ref{modernEP} above. Since there are $p-1$ distinct $p-1^{st}$ roots of unity, Dirichlet continued:
\begin{quote}
The equation just found represents $p-1$ different equations, which are obtained by replacing $\omega$ with its $p-1$ values. It is known that these $p-1$ different values can be represented as powers of one such $\Omega$, chosen appropriately, so that the values are then:
\[
\Omega^{0},\ \Omega^{1},\ \Omega^{2},\ \ldots,\ \Omega^{p-2}.
\]
In accordance with this representation, we will write the different values $L$ of the series or product as:
\[ 
L_{0},\ L_{1},\ L_{2},\ \ldots,\ L_{p-2},
\]
\ldots\footnote{``Die eben gefundene Gleichung repr\"{a}sentirt $p-1$ verschiedene Gleichungen, welche man erh\"{a}lt, wenn man f\"{u}r $\omega$ seine $p-1$ Werthe setzt.  Bekanntlich lassen sich diese $p-1$ verschiedenen Werthe durch die Potenzen von einem derselben $\Omega$ darstellen, wenn dieser geh\"{o}rig gew\"{a}hlt wird, und sind dann:\[
\Omega^{0},\ \Omega^{1},\ \Omega^{2},\ \ldots,\ \Omega^{p-2}.
\]
Wir werden, dieser Darstellung entsprechend, die verschiedenen Werthe $L$ der Reihe oder des Productes mit:
\[
L_{0},\ L_{1},\ L_{2},\ \ldots,\ L_{p-2}
\]
bezeichnen\ldots .'' \cite[p.~318]{Dirichlet37}} \cite[p.~3]{Dirichlet37}
\end{quote}
Notice that Dirichlet says that the Euler product formula ``represents $p-1$ different equations,'' rather than thinking of it as a single equation parametrized by $\omega$.

In the more general case where the common difference is some composite $k$, Dirichlet's procedure is completely analogous.  First, he demonstrated that the Euler product formula holds:
\begin{quote}
\begin{align}
 \prod \frac{1}{1-\theta^{\alpha}\ph^{\beta}\omega^{\gamma}\omega^{'\gamma^{'}}\ldots\frac{1}{q^{s}}}=\sum \theta^{\alpha}\ph^{\beta}\omega^{\gamma}\omega^{'\gamma^{'}}\ldots\frac{1}{n^{s}} = L, \label{dirichletEulergeneral}
\end{align}
where the multiplication sign ranges over all primes, with the exclusion of $2, p, p', \ldots$, and the summation ranges over all the positive integers that are not divisible by any of the primes $2, p, p', \ldots$. The system of indices $\alpha, \beta, \gamma, \gamma', \ldots$ on the left side corresponds to the number $q$, and on the right side to the number $n$. The general equation (\ref{dirichletEulergeneral}), in which the different roots $\theta, \ph, \omega, \omega', \ldots$ can be combined with one another arbitrarily, clearly contains $K$-many particular equations.\footnote{``\begin{align}
 \prod \frac{1}{1-\theta^{\alpha}\ph^{\beta}\omega^{\gamma}\omega^{'\gamma^{'}}\ldots\frac{1}{q^{s}}}=\sum \theta^{\alpha}\ph^{\beta}\omega^{\gamma}\omega^{'\gamma^{'}}\ldots\frac{1}{n^{s}} = L, \tag{\ref{dirichletEulergeneral}} \end{align}
 wo sich das Multiplicationszeichen auf die  ganze Reihe der Primzahlen, mit Ausschluss von 2, $p$, $p'$, \ldots, und das Summenzeichen auf alle positiven ganzen Zahlen, welche durch keine der Primzahlen 2, $p$, $p'$, \ldots theilbar sind, erstreckt.  Das System der Indices $\alpha, \beta, \gamma, \gamma', \ldots$ entspricht auf der ersten Seite der Zahl $q$, auf der zweiten Seite der Zahl $n$.  Die allgemeine Gleichung (\ref{dirichletEulergeneral}), in welcher die verschiedenen Wurzeln $\theta, \ph, \omega, \omega', \ldots$ auf irgend eine Weise mit einander combinirt werden k\"{o}nnen, enth\"{a}lt offenbar eine Anzahl $K$ besonderer Gleichungen.'' \cite[pp.~336--337]{Dirichlet37} We have replaced Dirichlet's equation number with our own.} \cite[p.~17]{Dirichlet37}
\end{quote}
Note, again, Dirichlet's characterization of the general equation as ``containing'' the particular instances.
Here, $K$ is what we have called $\ph(k)$, the cardinality of the group $(\ZZ / k\ZZ)^*$. Dirichlet went on to note that we can choose primitive roots of unity $\Theta, \Phi, \Omega, \Omega', \ldots$ so that all choices of $\theta, \ph, \omega, \omega', \ldots$ can be expressed as powers of these,
\[
\theta = \Theta^\mathfrak{a}, \ph = \Phi^\mathfrak{b}, \omega = \Omega^\mathfrak{c}, \omega' = \Omega^{\mathfrak{c}'}, \ldots,
\]
just as in the simpler case. He wrote that we can thus refer to the $L$-series in a ``convenient'' (``bequem'') way, as $L_{\mathfrak{a}, \mathfrak{b}, \mathfrak{c}, \mathfrak{c}', \ldots}$, where $\mathfrak{a}, \mathfrak{b}, \mathfrak{c}, \mathfrak{c}', \ldots$ are the exponents of the chosen primitive roots. Notice that the representations just described depend on fixed, but arbitrary, choices of the primitive roots of unity, as well as fixed but arbitrary generators of the cyclic groups. Modulo those choices, we have parameters $\mathfrak{a}, \mathfrak{b}, \mathfrak{c}, \mathfrak{c}', \ldots$ that vary to give us all the characters; and for each choice of $\mathfrak{a}, \mathfrak{b}, \mathfrak{c}, \mathfrak{c}', \ldots$ we have an explicit expression that tells us the value of the character at $n$. For Dirichlet, summing over characters therefore amounted to summing over all possible choices of this representing data.

In the special case of where the common difference is a prime, $p$, Dirichlet ran through calculations similar to those described in Section~\ref{contemporary:proofs:section} to obtain the following identity:
\begin{multline*}
\sum\frac{1}{q^{1 + \rho}} +\frac{1}{2}\sum\frac{1}{q^{2 +2 \rho}} + \frac{1}{3}\sum\frac{1}{q^{3 + 3 \rho}} + \ldots \\
= \frac{1}{p-1} (\log L_{0} + \Omega^{-\gamma_{m}} \log L_{1} + \Omega^{-2\gamma_{m}}\log L_{2} + \ldots + \Omega^{-(p-1)\gamma_{m}} \log L_{p-2}).
\end{multline*}
He has essentially arrived at equation (\ref{sumoverchar2}) in Section~\ref{contemporary:proofs:section}, which read as follows:
\begin{equation*}
\tag{\ref{sumoverchar2}}
\sum_{\chi \in \widehat{(\mathbb{Z} / k \mathbb{Z})^*}}\overline{\chi(m)}\log L(s, \chi) =
\ph(k) \sum_{q \equiv m \pmod{k}}\frac{1}{q^{s}} \ + \ O(1)
\end{equation*}
In the case at hand, $k$ is $p$, in which case $\ph(k)$, the cardinality of $(\mathbb{Z} / k \mathbb{Z})^*$, is equal to $p - 1$. To facilitate the comparison, switch the right-hand-side of (\ref{sumoverchar2}) with the left-hand-side (convention dictates that the term $O(1)$ stays on the right side), divide through by $\ph(k)$, and note all of the following. First, Dirichlet used $1 + \rho$ in place of $s$ to denote the quantity that approaches $1$ from above. Second, the sum
\[
\frac{1}{2}\sum\frac{1}{q^{2 +2 \rho}} + \frac{1}{3}\sum\frac{1}{q^{3 + 3 \rho}} + \ldots
\]
has been absorbed into the term $O(1)$ in equation (\ref{sumoverchar2}); indeed, Dirichlet's next move was to note that this sum is bounded by a constant. Finally, if $\chi_i(m)$ is equal to $\Omega^{i \gamma_m}$, then the complex conjugate $\overline{\chi_i(m)}$ is equal to $\Omega^{-i \gamma_m}$, so the expression $\Omega^{-i \gamma_m} \log L_i$ would be expressed in our notation as $\overline{\chi_i(m)} L(s, \chi_i)$.

When Dirichlet considered the more general case, he arrived at the analogous result:
\begin{multline*}
\sum\frac{1}{q^{1 + \rho}} + \frac{1}{2}\sum\frac{1}{q^{2 +2 \rho}} + \frac{1}{3}\sum\frac{1}{q^{3 + 3 \rho}} + \ldots \notag \\
= \frac{1}{K}\sum \Theta^{-\alpha_{m}\mathfrak{a}}\ \Phi^{-\beta_{m}\mathfrak{b}}\Omega^{-\gamma_{m}\mathfrak{c}}\Omega^{-\gamma'_{m}\mathfrak{c}'} \cdots \log L_{\mathfrak{a},\mathfrak{b},\mathfrak{c},\mathfrak{c}', \ldots.}
\end{multline*}
Here the summation on the right hand side of the equation is over the possible values of $\mathfrak{a}, \mathfrak{b}, \mathfrak{c}, \mathfrak{c'}, \ldots$. Once again, this translates to our equation (\ref{sumoverchar2}).

Finally, consider the key use of the first orthogonality relation given by Corollary~\ref{orthocorrol}. In the special case where the common difference is a prime, Dirichlet expressed this by saying that we have
\[
1+\Omega^{h\gamma - \gamma_{m}} + \Omega^{2(h\gamma - \gamma_{m})} + \ldots + \Omega^{(p-2)(h\gamma - \gamma_{m})}=0 
\]
except when $h\gamma - \gamma_{m} \equiv 0 \bmod{p-1}$, in which case the sum is equal to $p - 1$. In the general case, he considered the sum $\frac{1}{h}\sum W\frac{1}{q^{h + h\rho}}$,
\begin{quote}
\ldots where the symbol $\sum$ ranges over all primes $q$ and $W$ denotes the product of the sums taken over $\mathfrak{a}, \mathfrak{b}, \mathfrak{c}, \mathfrak{c'}, \ldots$ or respectively over 
\[                                                                                                       
\sum\Theta^{(h\alpha - \alpha_{m})\mathfrak{a}}, \sum\Phi^{(h\beta - \beta_{m})\mathfrak{b}}, \sum\Omega^{(h\gamma - \gamma_{m})\mathfrak{c}}, \sum\Omega'^{(h\gamma' - \gamma'_{m})\mathfrak{c'}}, \ldots
\]
Now one sees \ldots that the first of these sums is 2 or 0, corresponding to whether the congruence $h\alpha-\alpha_{m} \equiv 0 \pmod{2}$ or, equivalently, whether the congruence $q^{h} \equiv m \pmod{4}$ holds or not; that the second is $2^{\lambda-2}$ or 0 corresponding to whether the congruence $h\beta-\beta_{m} \equiv 0 \pmod{2^{\lambda-2}}$ or, equivalently, whether the congruence $q^{h} = \pm m \pmod{2^{\lambda}}$ holds or not; that the third is $(p-1)p^{\pi-1}$ or 0, corresponding to whether the congruence $h\gamma-\gamma_{m} \equiv 0 \pmod{(p-1)p^{\pi-1}}$ or, equivalently, whether the congruence $q^{h} \equiv m \pmod{p^{\pi}}$ holds or not, and so on; that therefore $W$ always vanishes except when the congruence $q^{h} \equiv m \pmod{k}$ holds, in which case $W=K$.\footnote{``\ldots wo sich das Zeichen $\sum$ auf die Primzahlen $q$ erstreckt und $W$ das Product der nach $\mathfrak{a}, \mathfrak{b}, \mathfrak{c}, \mathfrak{c'}, \ldots$ resp.\ zwischen den angegebenen Grenzen zu nehmenden Summen:
\[                                                                                                       
\sum\Theta^{(h\alpha - \alpha_{m})\mathfrak{a}}, \sum\Phi^{(h\beta - \beta_{m})\mathfrak{b}}, \sum\Omega^{(h\gamma - \gamma_{m})\mathfrak{c}}, \sum\Omega'^{(h\gamma' - \gamma'_{m})\mathfrak{c'}}, \ldots
\]
bedeutet. Nun erseiht man \ldots dass die erste dieser Summen 2 oder 0 ist, je nachdem die Congruenz $h\alpha-\alpha_{m} \equiv 0 \pmod{2}$, oder was dasselbe ist, die Conquenz $q^{h} \equiv m \pmod{4}$ stattfindet oder nicht stattfindet, das die zweite $2^{\lambda-2}$ oder 0 ist, je nachdem die Congruenz $h\beta-\beta_{m} \equiv 0 \pmod{2^{\lambda-2}}$, oder was dasselbe ist, die Conquenz  $q^{h} = \pm m \pmod{2^{\lambda}}$ oder was dasselbe ist, die Conquenz $q^{h} \equiv m \pmod{4}$ stattfindet oder nicht stattfindet, das die dritte $(p-1)p^{\pi-1}$ oder 0 ist, je nachdem die Congruenz $h\gamma-\gamma_{m} \equiv 0 \pmod{(p-1)p^{\pi-1}}$ oder was dasselbe ist, die Conquenz  $q^{h} \equiv m \pmod{p^{\pi}}$ stattfindet oder nicht stattfindet, u.~s.~w.~das also $W$ immer verschwindet, ausser wenn die Congruenz $q^{h} \equiv m \pmod{k}$ ist, in welchem Falle $W=K$ wird.''} \cite[p.~340]{Dirichlet37}
\end{quote}

Let us now summarize the differences between Dirichlet's presentation and contemporary ones. One salient difference is at the level of algebraic abstraction. The contemporary presentation developed the general notion of a character of an arbitrary finite abelian group, $G$, and showed that the characters themselves form a group, $\widehat G$, isomorphic to $G$ itself. Dirichlet, in contrast, focused on a very particular group, $(\ZZ/ k \ZZ)^*$, the multiplicative group of units modulo $k$. His presentation showed an intimate familiarity with the structure of that group, and the explicit mapping from that group to the group of characters. These details are inessential in the modern presentation.

But there is another abstraction that distinguishes modern presentations from Dirichlet's, namely, the willingness to treat characters as mathematical objects in their own right. This is, in part, facilitated by the algebraic abstraction: it is easier and more advantageous to treat characters as objects that transcend their representations when there are language and methods available that obviate the need to return to particular representations whenever there is real work to be done. But the dependence goes both ways: the abstract algebraic methods cannot even be invoked until one is willing to consider characters as objects that can bear properties and algebraic structure.

It is the absence of these two forms of abstraction that puts Dirichlet's presentation in stark contrast with the modern one. Dirichlet did not define the notion of character at all, let alone general operations on characters; nor did he identify any of their general properties. Rather, characters only come into play as symbolic expressions in the construction of the $L$-functions, and their properties are derived in an ad-hoc way as the proof proceeds. The characters are thus ``intensional'' entities: they are represented explicitly in terms of primitive roots of unity and generators of the corresponding residue group, and the only way of reasoning about them is in terms of this representing data.

In Dirichlet's presentation, there is no quantifying over characters. Instead, one quantifies over their representations in terms of primitive roots of unity. Similarly, there is no direct notion of summing over characters: Dirichlet was happy to sum over finite sets of natural numbers and tuples of natural numbers, but where we would sum over a finite set of characters, Dirichlet instead summed over representations in terms of such tuples. It is notable that he presented the Euler product formula in (\ref{dirichletEulersimple}) and (\ref{dirichletEulergeneral}) in terms of arbitrary primitive roots of unity, and then explained that these ``represent'' or ``contain'' more particular equations which can then be summed over to obtain the desired result.

Finally, it is worth drawing attention to one further consequence of the different treatment of characters in the two presentations, namely, the extent to which they make important dependences explicit. Because Dirichlet's expressions depend on so much representational data, Dirichlet often suppressed details, to keep the expressions from getting unwieldy. Thus, he wrote $\omega^\gamma$, suppressing the dependence of $\gamma$ on $n$, where we would write $\chi(n)$ for the character $\chi$ corresponding to $\omega$. This places a greater burden on the text of the proof, and the reader's memory, to keep track of the relevant dependences, and, for example, the ranges of a summation. Moreover, the modern notation $L(s,\chi)$ makes it easy to track the dependence on the character $\chi$, something that is lost in Dirichlet's presentation. Thus in the modern expression $\overline{\chi(m)} L(s, \chi)$ the role of $\chi$ is clear in both terms in the product; in Dirichlet's expression, $\Omega^{-i \gamma_m} \log L_i$ the connection is buried in the definition of $L_i$.

\section{Later presentations}
\label{evolution:section}

In this section, we will examine proofs and discussions of Dirichlet's theorem by Dedekind \cite{dirichlet:63b}, de la Vall\'{e}e-Poussin \cite{Poussin96,Poussin97}, Hadamard \cite{Hadamard96}, Kronecker \cite{kronecker:01}, and Landau \cite{landau1909,landau:27}, as well as proofs of related results by Dirichlet \cite{Dirichlet40,Dirichlet41} and Weber \cite{weber:82}. The proofs of the original theorem all share the same structure, so we will focus on describing the similarities and differences between the presentations of characters and $L$-functions, rather than describing each in full.

\subsection{Dirichlet}
\label{dirichlet:extensions:section}

Over the next few years, Dirichlet extended the methods used to prove Theorem~\ref{dirichlet:theorem} in two different directions, first, to consider primes represented by quadratic forms, and, second, to consider arithmetic progressions in the complex numbers. We will consider each of these extensions briefly.

Regarding the first extension, some background on quadratic forms will be helpful. Number theory has long been concerned with questions as to which numbers can be represented by a given algebraic expression in which the variables are taken to range over the integers. It is easy to characterize the numbers that can be represented by a linear form $a x + b$ in one variable, and quadratic reciprocity addresses the problem of which numbers can be represented by a quadratic form $a x^2 + b x + c$. When it comes to quadratic forms in two variables, things become more difficult. Fermat's famous theorem that the odd prime numbers that can be represented by the form $x^2 + y^2$ are exactly the ones that are congruent to one modulo four is considered a gem of number theory. Euler proved this and other of Fermat's claims in the eighteenth century, and Lagrange later extended the theory of binary quadratic forms considerably. A central contribution of Gauss' \emph{Disquitiones Arithmeticae} is a complete classification of binary quadratic forms, whose study constitutes the bulk of that work. 

Gauss showed that it suffices to characterize ``primitive'' quadratic forms $a x^2 + 2 b x y +c y^2$, where the second coefficient is even, and $a$, $b$, and $c$ have no factor in common. (Quadratic forms are also classified as ``indefinite,'' ``positive definite,'' and ``negative definite,'' and in the discussion that follows it should be assumed that we are fixing our attention on one of these fixed kinds.) He called the value $D = b^2 - a c$ the \emph{discriminant} of the form, and showed how to assign to each primitive quadratic form of discriminant $D$ a finite list of values of the form $\pm 1$ which, roughly speaking, characterizes its behavior. He called these values the \emph{characters} of the form, and took the ``genus'' of a form to be the collection of all the primitive forms with the same discriminant and character. As we will see in Secction~\ref{dedekind:section}, the use of name ``character'' for nonzero group homomorphisms to the complex numbers traces back to this terminology.

Dirichlet's Theorem~\ref{dirichlet:theorem} asserts that there are infinitely many primes represented by any \emph{linear} form $k x + m$ where $k$ and $m$ have no common factor. In 1840, Dirichlet considered the analogous question as to whether infinitely many primes can be represented by a primitive quadratic form $a x^2 + 2 b x y + c y^2$, and showed that this indeed holds in the special case where the discriminant $D = -p$ is a negative prime number, and $p$ is congruent to 3 mod 4. As he noted, the proof relies on methods that are ``essentially the same'' (``im Wesentlichen \ldots \"{u}bereinstimmend'') \cite[p.~98]{Dirichlet40} as those used in his proof of the theorem on arithmetic progressions, although with a few new modifications. The result was proved in full generality by Weber \cite{weber:82} in 1882, and de la Vall\'{e}e Poussin \cite{Poussin97} also presented a proof in 1897.

In 1841, Dirichlet turned his attention to arithmetic progressions in the Gaussian integers, and established that
\begin{quote}
\ldots the expression $kt+l$, in which $t$ denotes an indeterminate complex integer and $k, l$ denote given numbers of this kind with no common factor, always contains infinitely many prime numbers.''\footnote{``\ldots der Ausdruck $kt+l$, in welchem $t$ eine unbestimmte complexe ganze Zahl und $k, l$ gegebene solche Zahlen ohne gemeinschaftlichen Factor bezeichnen, immer unendlich viele Primzahlen enth\"{a}lt.''} \cite[p.~509]{Dirichlet41} 
\end{quote}
The proof begins with a discussion of modular arithmetic for the Gaussian integers, and defines the notion of primitive roots and indices for complex moduli \cite[p.~512--524]{Dirichlet41}. Dirichlet then considered the product 
\[
\Omega_{n}=\varphi^{\alpha_{n}}\varphi^{'\alpha'_{n}}\ldots\times\psi^{\beta_{n}}\chi^{\gamma_{n}}\psi^{'\beta'_{n}}\chi^{'\gamma'_{n}}\ldots\times\theta^{\delta_{n}}\eta^{\varepsilon_{n}}
\]
where $\varphi, \varphi', \psi, \psi', \chi, \chi', \theta, \eta$ are all roots of unity and $\alpha_{n}, \alpha'_{n}, \beta_{n}, \beta'_{n}, \chi_{n}, \chi'_{n}, \ldots, \allowbreak \delta_{n}, \varepsilon_{n}$ are indices with respect to $n$.  This product is thus analogous to the characters that appear in his 1837 proof.

However, in contrast to the presentation in 1837, in 1841 Dirichlet explicitly noted that the products $\Omega_{n}$ enjoy a number of ``important properties'' (``wichtige Eigenschaften'') \cite[p.~524]{Dirichlet41}.  The properties he highlights are the following, where $\psi(k)$ is a quantity calculated from the prime decomposition of $k$, involving the norms of the prime factors and their exponents \cite[524--525]{Dirichlet41}:
\begin{enumerate}
\item $\Omega_{nn'} = \Omega_{n}\Omega_{n'}.$
\item $\Omega_{n'}=\Omega_{n}$ whenever $n'\equiv n \pmod{k}$.
\item $\sum\Omega_{l} = 0$ or $\sum\Omega_{l} = \frac{1}{4}\psi(k)$ depending on whether there is at least one root among the roots in $\Omega_{l}$ that is different to 1, or whether they are all equal.
    \item $S\Omega_{n} = \frac{1}{4}\psi(k)$ or $S\Omega_{n}=0$ depending on whether $n\equiv 1 \pmod{k}$ or $n\not\equiv 1 \pmod{k}$, where the sign ``$S$'' indicates a sum over all combinations of the roots that can occur in $\Omega$.
\end{enumerate}
Note that the last two are the key orthogonality relations, Theorem~\ref{ortho}. Thus, in this presentation Dirichlet labeled the expressions that we now recognize as characters and isolated their properties. Moreover, he introduced a notation for summing over characters that is, interestingly, distinct from the notation for summation over natural numbers. The net effect is that the proof is considerably more modular than the proof of 1837.

\subsection{Dedekind}
\label{dedekind:section}

In 1863, Dedekind published a book, \emph{Vorlesungen \"{u}ber Zahlentheorie}, based on his notes taken from a course on number theory given by Dirichlet at G\"ottingen. After the lectures, he added nine ``supplements,'' or appendices, with material of his own. Supplement VI, in particular, contained a presentation of Dirichlet's theorem. Dedekind extended the supplements in the second edition, published in 1871, to include his theory of algebraic ideals. That theory was, in turn, revised and expanded in the third and fourth editions, which appeared in 1879 and 1894, respectively. The earlier supplements, however, were not changed after the first edition.

Dedekind began his presentation with an overview of main steps of the proof. In particular, he proved the Euler product formula, and obtained the series expansions for the $L$-functions and their logarithms. But the overview deals with a more general class of $L$-functions than we have considered so far:
\begin{quote}
The general proof of $\ldots$ [Dirichlet's theorem] is based on the consideration of a class of infinite series of the form
\[
L = \sum \psi(n),
\]
where $n$ runs through all positive integers and the real or complex function $\psi(n)$ satisfies the condition
\[
\psi(n)\psi(n')=\psi(nn'),
\]
\ldots [and] we always assume that $\psi(1)=1$.\footnote{``Der allgemeine Beweis dieses Satzes \ldots st\"{u}tzt sich auf die Betrachtung einer Classe von unendlichen Reihen von der Form \[
L = \sum \psi(n),
\]
wo der Buchstabe $n$ alle ganzen positiven Zahlen durchlaufen muss, und die reelle oder complexe Function $\psi(n)$ der Bedingung
\[
\psi(n)\psi(n')=\psi(nn')
\]
gen\"{u}gt \ldots so nehmen wir immer an, dass $\psi(1) = 1$ ist.''} \cite[\S 132]{dirichlet:63b}
\end{quote}

He then went on to focus his attention on Dirichlet characters. He used the same expressions as Dirichlet to represent the characters, but whereas Dirichlet expressed each root of unity $\omega$ involved in the construction as the power $\Omega^i$ of a single primitive root of unity, Dedekind did not bother with this step. He did not explicitly refer to the relevant expressions as ``characters,'' but he introduced the notation $\chi(n)$ to denote their values, and pointed out that $\chi$ is completely multiplicative:
\begin{quote}
The numerator [of $\psi(n)$] $\chi(n)=\theta^{\alpha}\eta^{\beta}\omega^{\gamma}\omega'^{\gamma'}\ldots$ has the characteristic property $\chi(n)\chi(n')=\chi(nn')\ldots$\footnote{``Der Z\"{a}hler $\chi(n) = \theta^{\alpha}\eta^{\beta}\omega^{\gamma}\omega'^{\gamma'}\ldots$ besitzt die charakteristischen Eigenschaften $\chi(n)\chi(n') = \chi(nn')$ \ldots''} \cite[\S 133, footnote]{dirichlet:63b}
\end{quote}

Otherwise, however, Dedekind did not make use of the $\chi$ notation in his proof. And, like Dirichlet, Dedekind showed a reluctance to quantify over characters directly. He introduced the $L$-functions as
\[
\prod\frac{1}{1-\psi(q)} = \sum\psi(n) = L,
\]
specifying that the product is to be taken over all primes not dividing the common difference, and the sum is to be taken over all natural numbers relatively prime to the common difference. Dedekind went on to note
\begin{quote}
\ldots that the series can exhibit quite different behavior, depending on the roots of unity $\theta, \eta, \omega, \omega', \ldots$ appearing in the expression for $\psi(n)$.  Since these roots can have $a, b, c, c', \ldots$ values, respectively, the form $L$ contains altogether
\[
abcc'\ldots = \ph(k)
\]
different particular series\ldots\footnote{``Wir bemerken zun\"{a}chst, dass diese Reihen je nach der Wahl der in dem Ausdrucke $\psi(n)$ vorkommenden Einheits-Wurzeln $\theta, \eta, \omega, \omega', \ldots$ ein ganz verschiedenes Verhalten zeigen; da diese Wurzeln resp. $a,b,c,c', \ldots$ verschiedene Werthe haben k\"{o}nnen, so sind in der Form $L$ im Ganzen
\[
abcc'\ldots = \ph(k)
\]
verschiedene besondere Reihen enthalten \ldots''} \cite[\S 133]{dirichlet:63b}
\end{quote}
Note that the use of the word ``contains'' echoes Dirichlet's language.

Dedekind then proceeded to divide the $L$-functions into classes.  Given that the characters are defined in terms of a sequence of roots of unity,  there are three distinct possibilities that can obtain for a given series $L$:
\begin{enumerate}
\item The roots of unity in the construction of the character occurring in $L$ are all 1.  There is only one such $L$-function, which is denoted $L_{1}$.
\item The roots of unity in the construction of the character occurring in $L$ are all real, i.e.\ are all $\pm 1$.  $L$-functions that fall into this category are written as $L_{2}$.
\item At least one of the roots of unity in the construction of the character occurring in $L$ is imaginary.  $L$-functions that fall into this category are written as $L_{3}$.
\end{enumerate}
Moreover, each $L$-function that falls into the third category has a conjugate: if $L_{3} = \sum \frac{\theta^{\alpha}\eta^{\beta}\omega^{\gamma}\omega'^{\gamma'}\ldots}{n^{s}}$, then there is a corresponding $L$-series in the same class, denoted by $L_{3}'$, such that $L_{3}' = \frac{\theta^{-\alpha}\eta^{-\beta}\omega^{-\gamma}\omega'^{-\gamma'}\ldots}{n^{s}}$.

Let us now see how Dedekind's notation plays out when it comes to summing over the characters.  When proving that the $L$-functions corresponding to a complex character have a finite non-zero value as $s$ tends to infinity, he obtained the following equation:
\begin{quote}
\begin{multline}
\label{Dedekindsummationcharacters1}
\ph(k)\left(\sum\frac{1}{q^s} + \frac{1}{2}\sum\frac{1}{q^{2s}} + \ldots + \frac{1}{\mu}\sum\frac{1}{q^{\mu s}} + \ldots\right) \\
 =\log L_{1} + \sum \log(L_{2}) + \sum \log(L_{3}L_{3}'),
\end{multline}
where, on the left hand side, the successive sums are over all the primes $q$ not dividing $k$ which satisfy the successive conditions $q \equiv 1 \pmod{k}$, $q^{2} \equiv 1 \pmod{k}$, etc.  On the right hand side the first sum is over all series $L_{2}$ of the second class, and the second sum is over all conjugate pairs $L_{3}L_{3}'$ of series of the third class.\footnote{``\[
\ph(k)\left(\sum\frac{1}{q^s} + \frac{1}{2}\sum\frac{1}{q^{2s}} + \ldots + \frac{1}{\mu}\sum\frac{1}{q^{\mu s}} + \ldots\right)
= \log L_{1} + \sum \log(L_{2}) + \sum \log(L_{3}L_{3}'),
\]
wo auf der linken Seite das erste, zweite Summenzeichen u.s.f. sich auf  alle die in $k$ nicht aufgehenden Primzahlen $q$ bezieht, welche resp.~den Bedingungen $q \equiv 1, q^{2} \equiv 1 \pmod{k}$ u.s.f. Gen\"{u}ge leisten; auf der rechten Seite bezieht sich das erste Summenzeichen auf alle Reihen $L_{2}$ der zweiten Classe, das zweite auf alle verschiedenen Paare $L_{3}L_{3}'$ conjugirter Reihen dritter Classe.'' We have added the equation numbers in this quotation and the next, for later reference.} \cite[\S 136]{dirichlet:63b}
\end{quote}
In a sense, Dedekind took the summations to range not over the characters or their representations, but over the \emph{$L$-functions} themselves. Later, when dealing with these sums, he came closer to Dirichlet's presentation, in that the indices of the summation range over the sequences of roots that determine the characters. However, at times he introduced convenient abbreviations. For example, given a particular collection of roots of unity $\theta, \eta, \omega, \omega', \ldots$, Dedekind denoted $\theta^{-\alpha_{1}}\eta^{-\beta_{1}}\omega^{-\gamma_{1}}\omega'^{-\gamma'_{1}}\ldots$ by $\chi$, where $\alpha_{1}, \beta_{1}, \gamma_{1}, \gamma'_{1}$ stand for the indices of the first term of the progression,  $m$. Note that while Dedekind used the notation $\chi(n)$ in a footnote, he explicitly called $\chi$, as defined here, a \emph{value}. (In fact, it plays the role of $\overline{\chi(m)}$ in the presentation of Section~\ref{dirichlet:characters:section}.) Moreover, when he used the symbol in a summation, Dedekind was explicit that the summation ranges not over $\chi$, but rather the roots of unity involved in the definition:
\begin{quote}
The summation of all products $\chi\log L$ therefore gives the result
\begin{align}
\label{Dedekindsummationcharacters2}
\ph(k)\left(\sum\frac{1}{q^{s}}+\frac{1}{2}\sum\frac{1}{q^{2s}} + \frac{1}{3}\sum\frac{1}{q^{3s}} + \ldots\right) = \sum \chi \log L,
\end{align}
where the successive sums on the left hand side are over all primes $q$ satisfying the successive conditions $q \equiv m, q^{2} \equiv m, q^{3} \equiv m \pmod{k}$ etc., while the sum on the right hand side is over all $\ph(k)$ different root systems $\theta, \eta, \omega, \omega', \ldots$\footnote{``Die Summation aller Producte $\chi\log L$ giebt daher das Resultat
\[
\ph(k)\left(\sum\frac{1}{q^{s}}+\frac{1}{2}\sum\frac{1}{q^{2s}} + \frac{1}{3}\sum\frac{1}{q^{3s}} + \ldots\right) = \sum \chi \log L,
\]
wo auf der linken Seite das erste, zweite, dritte Summenzeichen u.s.f. sich auf alle Primzahlen $q$ bezieht, welche resp. den Bedingungen $q \equiv m, q^{2} \equiv m, q^{3} \equiv m \pmod{k}$ u.s.f. gen\"{u}gen, w\"{a}hrend das Summenzeichen auf der rechten Seite sich auf die s\"{a}mmtlichen $\ph(k)$ verschiedenen Wurzel-Systeme $\theta, \eta, \omega, \omega', \ldots$''} \cite[\S 137]{dirichlet:63b}
\end{quote}

In many ways, Dedekind did not stray far from Dirichlet's presentation. For the most part, his treatment of the characters was intensional, in the sense that the arguments rely on the particular representations of the characters. In other words, operations on the characters are described in terms of their representations, rather than their values; and, like Dirichlet, he viewed summations as ranging over these representations in equations (\ref{Dedekindsummationcharacters1}) and (\ref{Dedekindsummationcharacters2}). Similarly, he classified the characters depending on the roots involved in their representation, rather than their values. Nor did he take the expression $L$, for the $L$-functions, to depend on the characters themselves. Rather, he characterized them as expressions that behave differently ``depending on the roots of unity'' appearing in their construction. The $L$ notation in particular does nothing to signal this dependence.

Nonetheless, he did take some key steps towards viewing the characters abstractly. To start with, Dedekind went out of his way to isolate the characters as independent of the $L$-series in which they appear, and, in particular, flagged them as entities satisfying certain key properties. And, at least at times, he characterized summations as ranging over the $L$-functions themselves, hinting at a new level of abstraction.

The use of the term ``character'' in the modern sense can be traced to the long Supplement XI in the 1879 edition of the Dirichlet-Dedekind \emph{Vorlesungen}, which contains a presentation of Dedekind's theory of ideals in algebraic number fields. As noted in Section~\ref{dirichlet:extensions:section}, Gauss assigned ``characters'' to quadratic forms in order classify their behavior. Already in the first edition of the \emph{Vorlesungen}, Dedekind used Gauss' terminology in Supplement IV, ``Genera of quadratic forms.'' In 1879, however, Dedekind went further by showing that Gauss' classification could be understood in terms of his theory of ideals. Specifically, he showed that there is a correspondence between genera of quadratic forms and equivalence classes of ideals in an associated quadratic extension of the rationals. Moreover, there is a group structure on the ideals, and Gauss' characters correspond to characters on that group, in the modern sense. This explains the quotation in Section~\ref{functions:section}: Dedekind had simply adopted Gauss' terminology to characterize the homomorphisms from groups of ideals to the (nonzero) complex numbers more generally.

As we noted in Section~\ref{generalization:section}, in 1882, Weber \cite{weber:82}, building on Dedekind's work, extended the notion of a ``character'' to arbitrary finite abelian groups. But although Dedekind revised his theory of ideals substantially in the third and fourth editions of the \emph{Vorlesungen}, he did not revise any of the supplements that appeared in the first edition. In particular, he did not take advantage of the opportunity to go back to introduce the modern notion of a character in his presentation of Dirichlet's proof, although he clearly could have done so in the later editions.

\subsection{Weber}
\label{weber:section}

Weber's 1882 paper \cite{weber:82} is devoted to a proof of the fact that any primitive quadratic form represents infinitely many primes, the result that Dirichlet had proved in 1840 in a special case. That paper, however, was the first to present the modern notation of a character on a finite abelian group in full generality. Weber defined such a group as consisting of $h$ elements ``of any type'' (``irgend welcher Art'') satisfying the usual group laws. He then stated the now-familiar structure theorem:
\begin{quote}
In a finite abelian group of order $h$ one can always choose elements $\Theta_{1}, \Theta_{2}, \ldots, \Theta_{\nu}$ of order $n_{1}, n_{2}, \ldots, n_{\nu}$ so that each element $\Theta$ of $G$ can be expressed uniquely by the form $$\Theta_{1}^{s_{1}}\Theta_{2}^{s_{2}}\ldots \Theta_{\nu}^{s_{\nu}},$$ where $s_{1}, s_{2}, \ldots, s_{\nu}$ are chosen from complete residue systems with respect to the moduli $n_{1}, n_{2}, \ldots, n_{\nu}$.\footnote{``In einer Abel'schen Grupper $G$ von Grade $h$ kann man stets die Elemente $\Theta_{1}, \Theta_{2}, \ldots, \Theta_{\nu}$ von den Graden $n_{1}, n_{2}, \ldots, n_{\nu}$ so ausw\"{a}hlen, dass in der Form $$\Theta_{1}^{s_{1}}\Theta_{2}^{s_{2}}\ldots \Theta_{\nu}^{s_{\nu}}$$ jedes Element $\Theta$ von $G$ und jedes nur einmal enthalten ist, wenn $s_{1}, s_{2}, \ldots, s_{\nu}$ je einem vollst\"{a}ndigen Restsystem nach den Moduln  $n_{1}, n_{2}, \ldots, n_{\nu}$ entnommen werden.''} \cite[pp.~306--307]{weber:82}
\end{quote}

Weber called the sequence of values $\Theta_{1}, \Theta_{2}, \ldots, \Theta_{\nu}$ a \emph{basis} of the group, and used them to define the characters as follows:
\begin{quote}
If we assign to the $\nu$ elements $\Theta_{1}, \Theta_{2}, \ldots, \Theta_{\nu}$ of such a basis $\nu$ roots of unity $\omega_{1}, \omega_{2}, \ldots, \omega_{\nu}$, of order $n_{1}, n_{2}, \ldots, n_{\nu}$ respectively, then each element $\Theta =\Theta_{1}^{s_{1}}\Theta_{2}^{s_{2}}\ldots \Theta_{\nu}^{s_{\nu}}$ of the group also corresponds to a particular $h^{\mathit{th}}$ root of unity $$\omega =\omega_{1}^{s_{1}}\omega_{2}^{s_{2}}\ldots \omega_{\nu}^{s_{\nu}}.$$  We denote this root of unity by $\chi(\Theta)$ and call it the 
\emph{character} of the element.\footnote{``Ordnet man den $\nu$ Elementen $\Theta_{1}, \Theta_{2}, \ldots, \Theta_{\nu}$ einer solchen Basis $\nu$ Einheitswurzeln $\omega_{1}, \omega_{2}, \ldots, \omega_{\nu}$ von den Graden $n_{1}, n_{2}, \ldots, n_{\nu}$, so entspricht auch jedem Element $\Theta = \Theta_{1}^{s_{1}}\Theta_{2}^{s_{2}}\ldots \Theta_{\nu}^{s_{\nu}}$ der Gruppe eine bestimmte $h^{\mathit{te}}$ Einheitswurzel $\omega$ nach der Vorschrift $$\omega =\omega_{1}^{s_{1}}\omega_{2}^{s_{2}}\ldots \omega_{\nu}^{s_{\nu}}.$$  Wir bezeichnen diese Einheitswurzel mit $\chi(\Theta)$, und nennen dieselbe den \emph{Charakter} des Elementes $\Theta$.''}  \cite[p.~307]{weber:82}
\end{quote}
Weber went on to point out that this gives rise to $h$ distinct characters, and that each such character $\chi$ satisfies the equation
\[
 \chi(\Theta) \chi(\Theta') = \chi(\Theta \Theta').
\]
Moreover, this last condition provides an exact characterization:
\begin{quote}
If, conversely, $\chi(\Theta)$ is a uniquely determined function of $\Theta$ which satisfies [the equation above], then it is necessarily contained among these $h$ characters.\footnote{Ist umgekehrt $\chi(\Theta)$ eine durch das Element $\Theta$ eindeutig bestimmte Function, welche der Bedingung \ldots gen\"ugt, so ist dieselbe nothwending unter diesen $h$ Charakteren enthalten.}  [\emph{ibid.}]
\end{quote}
Thus Weber's paper provides us not only with the modern notion of a character defined on an arbitrary finite abelian group, but also with the modern understanding that such characters constitute an instance of the function concept.

There are some small differences between Weber's presentation and ours. In his presentation of the second orthogonality relation, Weber used ellipses rather than summation notation to sum over the characters:
\begin{quote}
For each element $\Theta$:
\[
\chi_{1}(\Theta) + \chi_{2}(\Theta) + \ldots + \chi_{h}(\Theta) = 0,
\]
except for the identity element $\Theta_{0}$, for which we have
\[
\chi_{1}(\Theta_{0}) + \chi_{2}(\Theta_{0}) + \ldots + \chi_{h}(\Theta_{0}) = h.\footnote{ ``F\"{u}r jedes Element $\Theta$ ist:
\[
\chi_{1}(\Theta) + \chi_{2}(\Theta) + \ldots + \chi_{h}(\Theta) = 0,
\]
ausgenommen f\"{u}r das Hauptelement $\Theta_{0}$, f\"{u}r welches
\[
\chi_{1}(\Theta_{0}) + \chi_{2}(\Theta_{0}) + \ldots + \chi_{h}(\Theta_{0}) = h.
\]}
\]
\cite[p.~308]{weber:82}
\end{quote}
This implicitly assumes that the $h$ characters are enumerated $\chi_1, \ldots, \chi_h$, though the enumeration is arbitrary. Moreover, Weber provided the explicit construction of the set of characters before the extensional characterization, whereas modern presentations usually provide the extensional characterization first. But these points are minor, and, otherwise, his presentation is little different from the one in Section~\ref{group:character:section}.\footnote{ Mackey's historical survey \cite{mackey:80} of the history of harmonic analysis include a very helpful and informative overview of the history of character-theoretic ideas in number theory.} Even his fairly modern presentation harks back to an intensional view, however. After introducing analogues of the Euler product formula and the result of taking logarithms of both sides, he emphasized, echoing the language used by Dirichlet and Dedekind:
\begin{quote}
Each of the formulas \ldots represents $h$ different formulas, corresponding to the $h$ different characters $\chi_1, \chi_2, \ldots, \chi_h$.\footnote{``Jede der Formeln \ldots repr\"asentirt $h$ verschiedene Formeln, entsprechend den $h$ verschiedenen Charakteren $\chi_1, \chi_2, \ldots, \chi_h$.''} \cite[p.~315]{weber:82}
\end{quote}
He then characterized the analogue of equation (\ref{sumoverchar2}) in Section~\ref{contemporary:proofs:section} as the result of ``the addition of all the formulas'' (``die Addition der s\"ammtlichen Formeln''), rather than the addition of the corresponding values.

\subsection{De la Vall\'{e}e-Poussin}

More than two decades later, Charles Jean de la Vall\'{e}e-Poussin gave another presentation of Dirichlet's theorem in a 1895/6 paper entitled ``D\'{e}monstration simplifi\'{e}e du Th\'{e}or\`{e}me de Dirichlet sur la progression arithm\'{e}tique,'' and again in sections from his 1897 book \emph{Recherches analytiques sur la th{\'e}orie des nombres premiers}. He introduced the characters in much the same way that Dirichlet and Dedekind did, namely, via an explicit construction in terms of primitive roots and roots of unity. Like Dedekind and Weber, he used the symbol $\chi$ for the characters. Like Weber, he distinguished between the $\ph(M)$ characters modulo $M$ using subscripts, writing $\chi_{1}, \chi_{2}, \ldots, \chi_{\ph(M)}$.

De la Vall\'{e}e-Poussin went on to subject the characters to a thorough study. Whereas Dirichlet and Dedekind divided the $L$-functions (or, perhaps more precisely, the corresponding expressions) into three categories based on their representations, de la Vall\'ee Poussin based the categorization on the characters themselves. The following description applies to the simplest case, where the modulus $M$ is a prime:
\begin{quotation}
\noindent One calls the character that corresponds to the root $+1$ the \emph{principal character}; it is equal to unity for all numbers $n$. Apart from the principal character, there is only one which is real for all numbers $n$: it corresponds to the root ($-1$) and is equal to $\pm 1$ depending on the number $n$.

We call all the other characters \emph{imaginary characters}, though they may have a real value for some particular numbers. Their modulus is always equal to unity.\footnote{``On appelle \emph{caract\`{e}re principal} celui qui correspond \`{a} la racine $+1$; il est \'{e}gal \`{a} l'unit\'{e} pour tous les nombres $n$.  En dehors du caract\`{e}re principal, il n'y en a qu'un seul qui soit r\'{e}el pour tous les nombres $n$: il correspond \`{a} la racine ($-1$) et est \'{e}gal \`{a} $\pm 1$ suivant le nombre $n$.

Nous donnerons \`{a} tous les autres caract\`{e}res le nom de \emph{caract\`{e}re imaginaires,} quoiqu'ils puissent avoir une valeur r\'{e}elle pour certains nombres particuliers.  Leur module est toujours \'{e}gal \`{a} l'unit\'{e}.''} \cite[p.~19]{Poussin97}
\end{quotation}
Thus de la Vall\'{e}e Poussin divided the characters into three classes:
\begin{enumerate}
\item the class consisting solely of the principal character;
\item the class consisting of all real characters (in this case, there is only one character which is real for all numbers, corresponding to the root $-1$); and
    \item The class consisting of all other characters, called the imaginary characters.
\end{enumerate}
Notice that his categorization includes both intensional and extensional characterizations; that is, he characterized the characters in each class both in terms of the values they take, and the roots involved in their construction.

In his study of the characters in his 1896 paper, de la Vall\'{e}e Poussin listed a number of ``very important relations'' (``relations tr\`{e}s importantes''). The first is that $\chi(n)\chi(n')=\chi(nn')$ for every character $\chi$, and $n$ and $n'$. The second is that for a given character $\chi$ modulo $M$, and any $n$ and $n'$, if $n \equiv n' \pmod{M}$ then $\chi(n) = \chi(n')$. The third and fourth are the first and second orthogonality relations, respectively. In order to present the second orthogonality relation, de la Vall\'{e}e Poussin introduced a new symbol, $S$, to denote summation over the characters, just as Dirichlet did in 1841:\footnote{``Consid\'{e}rons \ldots la somme \'{e}tendue \`{a} tous les caract\`{e}res, c'est-\`{a}-dire \`{a} tous les syst\`{e}mes de racines
\[
S_{\chi}\chi(n)=S_{\omega}\omega_{1}^{\nu_{1}}\omega_{2}^{\nu_{2}} \ldots
\]
\ldots\ \emph{Pour tout nombre $n$, la somme \'{e}tendue \`{a} la totalit\'{e} des caract\`{e}res
\[
S_{\chi}\chi(n)=0,
\]
\`{a} la seule exception pr\`{e}s du cas o\`{u} $$n \equiv 1 \pmod{M},$$ car alors tous les indicateurs sont nuls et l'on a
\[
S_{\chi}\chi(n) = \ph(M).\mbox{''}
\]}}
\begin{quote}
Consider \ldots the sum extending over all the characters, that is to say over all the systems of roots
\[
S_{\chi}\chi(n)=S_{\omega}\omega_{1}^{\nu_{1}}\omega_{2}^{\nu_{2}} \ldots
\]
\ldots\ \emph{For every number $n$, the sum extending over the totality of characters satisfies
\[
S_{\chi}\chi(n)=0,
\]
the only exception being the case where
\[
n\equiv 1 \pmod{M},
\]
because then all the indices are zero and one has
\[
S_{\chi}\chi(n)=\ph(M).
\]}
\cite[pp.~14--15]{Poussin96}
\end{quote}
Although his notation is different from that of Dirichlet's original notation and Dedekind's, his proof of Theorem~\ref{dirichlet:theorem}, like theirs, relies on an explicit calculation based on the construction of the characters, in contrast to the modern proof sketched in Section~\ref{group:character:section}.

When introducing the $L$-functions in his 1896 paper, de la Vall\'{e}e-Poussin simply described them in terms of their equivalent expressions as a sum and product. In 1897, however, he adopted a functional notation:
\begin{quote}
We define the function $Z(s, \chi \bmod{M})$, for $\mathcal{R}(s)>1$, by the absolutely convergent expressions:
\[
Z(s, \chi) = \mathop{\displaystyle{\sum}'}_{n=1}^{\infty} \frac{\chi(n)}{n^{s}}= \prod\left(1-\frac{\chi(q)}{q^{s}}\right)^{-1}
\]
where $n$ designates successively all the integers prime to $M$ and $q$ all the prime numbers not dividing $M$.\footnote{``Nous d\'{e}finirons la fonction $Z(s, \chi \bmod{M})$, pour $\mathcal{R}(s)>1$, par les expressions absolument convergentes
\[
Z(s, \chi) = \mathop{\displaystyle{\sum}'}_{n=1}^{\infty} \frac{\chi(n)}{n^{s}}= \prod\left(1-\frac{\chi(q)}{q^{s}}\right)^{-1}
\]
o\`{u} $n$ d\'{e}signe successivement tous les nombres entiers premiers \`{a} $M$ et $q$ tous les nombres premiers qui ne divisent pas $M.$''} \cite[p.~56]{Poussin97}
\end{quote}
Aside from the choice of the letter $Z$ instead of the letter $L$, we have finally arrived at the modern notation. But, like Dirichlet and Dedekind, he was reticent to quantify over the characters, using language the is eerily reminiscent of theirs.  For example, in his 1897 work, he wrote:
\begin{quote}
\ldots one finds the \emph{fundamental equation}
\[
(E)\ldots -\lim_{s=1}(s-1)\frac{Z'(s, \chi)}{Z(s, \chi)}=\lim_{s=1}(s-1)\sum_{q}\chi(q)\frac{lq}{q^{s}},
\]
and this equation (E) represents in reality $\ph(M)$ distinct ones, which result from exchanging the characters amongst themselves.\footnote{``\ldots on trouve l'\emph{\'{e}quation fondamentale} 
\[
(E)\ldots -\lim_{s=1}(s-1)\frac{Z'(s, \chi)}{Z(s, \chi)}=\lim_{s=1}(s-1)\sum_{q}\chi(q)\frac{lq}{q^{s}},
\]
et cette \'{e}quation (E) en repr\'{e}sente en r\'{e}alit\'{e} $\ph(M)$ distinctes par l'\'{e}change des caract\`{e}res entre eux.''} \cite[p.~65]{Poussin97}.
\end{quote}

Nonetheless, there are a number of important respects in which de la Vall\'{e}e Poussin's presentation is close to the modern one. To start with, characters are treated as objects of study in their own right, bearing their own properties and relations. Moreover, they are classified extensionally, although this classification is, at the same time, related to properties of their representations. They are also represented notationally as arguments to the $L$-functions. Finally, de la Vall\'{e}e-Poussin used a summation notation with an index ranging over the characters, not their representing data. The fact that he went out of his way to use one symbol, $S$, for summation of characters, and the usual $\Sigma$ symbol for summation over sets of natural numbers hints that he does not conceive of these as precisely the same mathematical operation, although they have similar properties.

\subsection{Hadamard}

For each real number $x$, let $\pi(x)$ denote the number of primes less than $x$. Around the turn of the century, both Legendre and Gauss conjectured that $\pi(x)$ is asymptotic to $x / \ln(x)$, in the sense that their ratio, $\pi(x) \ln x / x$, approaches $1$ as $x$ approaches infinity. This fact, now known as the ``prime number theorem'' was finally proved by both de la Vall\'{e}e-Poussin and Jacques Hadamard, working independently, in 1896. De la Vall\'{e}e-Poussin's proof was first published in the \emph{Annales de la Soci\'et\'e Scientifique de Bruxelles}, but also appeared in the 1897 book we discussed in the last section. Both Hadamard and de la Vall\'ee-Poussin also obtained a generalization of the prime number theorem to arithmetic progressions, which says that if $m$ is relatively prime to $k$, the number of primes less than $x$ that are congruent to $m$ modulo $k$ is asymptotic to $(1 / \ph(k)) \cdot (x / \ln x)$. This is, of course, a generalization of Dirichlet's theorem, since it implies that there are infinitely many primes congruent to $m$ modulo $k$.

Hadamard's 1896 paper was titled ``Sur la distribution des z\'eros de la fonction $\zeta(s)$ et ses cons\'equences arithm\'etiques,'' and provided a number of results concerning $L$-functions and characters. He introduced characters in the same way as Dirichlet, Dedekind and de la Vall\'{e}e-Poussin, namely, by a construction in terms of roots of unity and primitive elements.  Like Weber and de la Vall\'{e}e-Poussin, he distinguished between the characters modulo $k$ notationally by making use of subscripts, writing $\psi_{v}(n)$ where $v$ runs from 1 to $\ph(k)$. He defined the $L$-functions as follows:
\[
L_{v}(s)=\sum_{n=1}^{\infty}\frac{\psi_{v}(n)}{n^{s}}.
\]
The notation $L_{v}(x)$, in contrast to de la Vall\'{e}e-Poussin's $Z(s,\chi)$, does not indicate an explicit dependence on the character, though the subscript provides a link. As did his predecessors, Hadamard classified the $L$-functions into three different classes, but his classification was intensional, like Dirichlet's and Dedekind's, referring to the roots of unity used in the construction of the characters.

When it came to summing over the characters, Hadamard, unlike Dirichlet, Dedekind, or de la Vall\'{e}e-Poussin, let the index of the summation range over the \emph{subscripts} described above.
\begin{quote}
The fundamental equation that Dirichlet used in the demonstration of his theorem is
\begin{equation*}
\sum_{v}\frac{\log L_{v}(s)}{\psi_{v}(m)}=\ph(k)\left(\sum\frac{1}{q^{s}} + \frac{1}{2}\mathop{\displaystyle{\sum}'} \frac{1}{q^{2s}} + \frac{1}{3}\mathop{\displaystyle{\sum}''}\frac{1}{q^{3s}} + \ldots \right),
\end{equation*}
where $m$ is some integer prime to $k$, and where, among the signs $\sum, \sum', \sum'', \ldots$, the first ranges over the prime numbers $q$ such that $q \equiv m \pmod{k}$, the second ranges over the primes numbers $q$ such that $q^{2} \equiv m \pmod{k}$, etc.\footnote{``L'\'{e}quation fondamental utilis\'{e}e par Dirichlet pour la d\'{e}monstration de son th\'{e}or\`{e}me, est
\[
\sum_{v}\frac{\log L_{v}(s)}{\psi_{v}(m)}=\ph(k)\left(\sum\frac{1}{q^{s}} + \frac{1}{2}\mathop{\displaystyle{\sum}'} \frac{1}{q^{2s}} + \frac{1}{3}\mathop{\displaystyle{\sum}''}\frac{1}{q^{3s}} + \ldots \right)
\]
o\`{u} $m$ est un entier quelconque premier avec $k$ et o\`{u} les signes $\sum, \sum', \sum'', \ldots$ s'\'{e}tendent, le premier aux nombres premiers $q$ tels que $q \equiv m \pmod{k}$, le second aux nombres premiers $q$ tels que $q^{2} \equiv m \pmod{k}$, etc.''} \cite[p.~209]{Hadamard96}
\end{quote}
Once again, this is comparable to the modern formulation in equation (\ref{sumoverchar2}).

Hadamard's presentation falls somewhere between those of Dirichlet and Dedekind, on the one hand, and that of de la Vall\'ee-Poussin, on the other. His treatment of characters was, for the most part, intensional. But, in contrast to Dirichlet and Dedekind, he left representational data out of the key summations, thereby eliminating unnecessary clutter. However, unlike de la Vall\'{e}e-Poussin, he did not go so far as to characterize these as summations over the characters themselves. Rather, he introduced natural number indices, $v$, labeling the characters, and then took the variable of summation to range over those.

\subsection{Kronecker}

Between 1863 and 1891, Leopold Kronecker lectured at the University of Berlin on a range of subjects, including number theory, algebra, and the theory of determinants. After Kronecker's death, his student, Kurt Hensel, edited the five volumes of his collected works, which were published between 1895 and 1930. Hensel also took it upon himself to work Kronecker's copious lecture notes and course material into two textbooks, \emph{Vorlesungen \"uber Zahlentheorie} and \emph{Vorlesungen \"uber die Theorie der Determinanten}, which he published in 1901 and 1903, respectively. The first of these closes with a proof of Dirichlet's theorem, which we will discuss here.

Kronecker, in fact, wrote his doctoral dissertation on algebraic number theory under Dirichlet's supervision, completing it in 1845. Kronecker insisted that mathematics should maintain a clear focus on symbolic representations and algorithms, a commitment that is evident throughout his work. Avoiding talk of ``arbitrary'' functions, real numbers, and so on, Kronecker focused instead on the construction of algebraic systems and explicit algorithms for calculating with these algebraic representations. For example, his article ``Grundz\"uge einer arithmetischen Theorie der algebraischen Gr\"ossen'' \cite{kronecker:82} provides means of carrying out operations on systems of algebraic integers in finite extensions of the rationals. Similarly, ``Ein Fundamentalsatz der allgemeinen Arithmetik'' \cite{kronecker:87b} provides an explicit construction of a splitting field for any polynomial with integer coefficients, and he viewed this as filling a gap in Galois' work.\footnote{See Edwards \cite{edwards:80,edwards:89,edwards:07,edwards:09} for a discussion of these works, and Kronecker's mathematics more generally.}

Kronecker's approach to number theory had a similar orientation. As Hensel put it in the introduction to \emph{Vorlesungen \"uber Zahlentheorie}:
\begin{quote}
He believed that one can and must in this domain formulate each definition in such a way that its applicability to a given quantity can be assessed by means of a finite number of tests. Likewise, an existence proof for a quantity is to be regarded as entirely rigorous only if it contains a method by which that quantity can actually be found. Kronecker was far from being one to completely reject a definition or proof that does not meet these highest requirements, but he believed that there was then something missing, and he held that completing it along these lines is an important task, by which our knowledge is furthered in an important sense.\footnote{``Er meinte, man k\"onne und man m\"usse in diesem Gebiete eine jede Definition so fassen, da{\ss} durch eine endliche Anzahl von Versuchen gepr\"uft werden kann, ob sie auf eine vorgelegte Gr\"o{\ss}e anwendbar ist oder nicht. Ebenso w\"are ein Existenzbeweis f\"ur eine Gr\"o{\ss}e erst dann als v\"ollig streng anzushen, wenn er zugleich ein Method enthalte, durch welche die Gr\"o{\ss}e, deren Existenz bewiesen werde, auch wirklich gefunden werden kann. Kronecker war weit davon entfernt, eine Definition oder einen Beweis vollst\"andig zu verwerfen, der jenen h\"ochsten Anforderungen nicht entsprach, aber er glaubte, da{\ss} dann eben noch etwas fehle, und er hielt eine Erg\"anzung nach dieser Richtung hin f\"ur eine wichtige Aufgabe, durch die unsere Erkenntnis in einem wesentlichen Punkte erweitert w\"urde.'' We have modified and extended a translation due to Stein \cite[p.~250]{stein:88}.} \cite[p.~vi]{kronecker:01}
\end{quote}

Kronecker's proof of Dirichlet's theorem is a focal point of the book, in that it brings together methods and ideas developed throughout the entire work. Here is Hensel's characterization:
\begin{quote}
 [The book] closes with the proof of the famous theorem that any arithmetic sequence, whose first term and common difference are relatively prime, contains infinitely many prime numbers. But Kronecker completed Dirichlet's proof of this theorem in a significant sense, in that he proved that one can determine, for an arbitrarily large number $\mu$, a larger number $\bar \mu$, so that in the interval $(\mu \cdots \bar \mu)$ one is sure to find a prime number of the required form. This nice supplement to that famous proof is a fruit of the higher demands, mentioned above, that Kronecker placed on arithmetic proofs. And here it seems, in fact, that with this improvement of Dirichlet's theorem, nothing by way of simplicity or transparency has been lost.\footnote{``\ldots schlie{\ss}t mit dem Beweise des ber\"uhmten Satzes, da{\ss} jede arithmetische Reihe, deren Anfangsglied und Differenz teilerfremd sind, unendlich viele Primzahlen enth\"alt; aber Kronecker vorvollst\"andigt den Dirichletschen Beweis dieses Satzes in einem wesentlichen Punkte, indem er nachweist, da{\ss} man f\"ur jede beliebige gro{\ss} anzunehmende Zahl $\mu$ eine gr\"o{\ss}ere Zahl $\bar \mu$ so bestimmen kann, da{\ss} in dem Intervalle $(\mu \cdots \bar \mu)$ sich sicher eine Primzahl der verlangten Form befindet. Dies sch\"one Erg\"anzung jenes ber\"uhmten Beweises ist eine Frucht der oben erw\"ahnten h\"oheren Forderungen, welche Kronecker an arithmetische Beweise stellte, und hier scheint es in der That, da{\ss} durch diese Verbesserung der Dirichletsche Beweis nichts an Einfachheit und Durchsichtigkeit verloren hat.''} \cite[p.~viii]{kronecker:01}
\end{quote}

The interaction between analytic and number-theoretic methods is a central theme of the \emph{Vorlesungen}, which opens with a discussion of Gauss' distinction between the fields of number theory, algebra, and analysis, and argues that they cannot be cleanly separated. For example, Kronecker praised Leibniz' characterization of $\pi$ in terms of the series
\[
 \frac{\pi}{4} = 1 - \frac{1}{3} +\frac{1}{5} - \frac{1}{7} + \cdots = \sum_{n = 0}^\infty \frac{(-1)^n}{2n + 1}
\]
as providing a definition of $\pi$ of a ``fully number-theoretic character'' (``durchaus zalhlentheoretischem Charakter'').
\begin{quote}
 What these examples teach us holds, more generally, for all the definitions of analysis. These always lead back to the integers and their properties, and from that entire branch of mathematics, only the concept of a limit has so far remained foreign. Arithmetic cannot be separated from analysis, which has freed itself from its original source, geometry, and has developed independently, on free soil. Even less so, as Dirichlet has succeeded in obtaining the deepest and most beautiful results in arithmetic by combining the methods of the two disciplines.\footnote{``Was uns diese Beispiele lehren, ist nun ma{\ss}gebend f\"ur alle Definitionen der Analysis \"uberhaupt. Dieselben f\"uhren stets auf die ganzen Zahlen und ihre Eigenschaften zur\"uck, und es ist von dem ganzen Gebiete des letzgenannten Zweiges der Mathematik der einzige Begriff des limes oder der Grenze der Zahlentheorie bisher fremd geblieben. Gegen die Anaylsis also, die sich von ihrer urspr\"unglichen Quelle, der Geometrie, befreit und auf freiem Boden selbst\"andig entwickelt hat, kann die Arithmetik nicht abgegrenzt werden, um so weniger, als es \emph{Dirichlet} gelungen ist, grade die sch\"onsten und tiefliegenden arithmetischen Resultate durch die Verbindung der Methoden beider Disciplinen zu erzielen.''} \cite[pp.~4--5]{kronecker:01}
\end{quote}

This interplay comes to the fore in the proof of Dirichlet's theorem. The text reports that Kronecker's version of the proof, in the case where the common difference is prime, was worked out in the lectures he gave in the winter semester of 1875/1876, whereas the general case was presented in the winter semester of 1886/1887 \cite[p.~442]{kronecker:01}. Fixing a modulus $m$ and an $r$ relatively prime to $m$, recall that Kronecker aimed to provide, for a given number $\mu$, an explicit upper bound $\bar \mu$ on the integers one has to consider to find a prime number greater than $\mu$ and congruent to $r$ modulo $m$. When it came to the analytic part of the proof, Kronecker explained that obtaining the desired bound is reduced to obtaining a positive lower bound on a certain analytic series that arises in the proof. He wrote:
\begin{quote}
For the ambiguous [real] characters, Dirichlet's proof meets this requirement. But his methods are not sufficient to do the same for the series corresponding to the complex characters.\footnote{``F\"ur die ambigen Charaktere erf\"ullt der Beweis von Dirichlet auch diese Forderung, dagegen reichen seine Methoden nicht aus, um dasselbe auch f\"ur die Reihen zu leisten, welche den complexen Charakteren entsprechen.''} \cite[p.~481]{kronecker:01}
\end{quote}
He went on:
\begin{quote}
Generally speaking, this is a special case of the problem of finding, given a well-defined nonzero number, a bound, above which it necessarily lies. This is not as easy as it seems at first glance; indeed, in some circumstances, the problem can count among the thorniest questions known to science.\footnote{``\"Uberhaupt ist das hier in einen speziellen Falle sich darbietende Problem, f\"ur eine von Null verschiedene wohldefinierte Zahlgr\"o{\ss}e eine Grenze zu finden, \emph{\"uber} der sie notwendig liegen mu{\ss}, nicht so einfach, als es auf den ersten Blick erscheint, vielmehr kann diese Aufgabe unter Umst\"anden eine der heikelsten Fragen sein, die die Wissenschaft kennt.''} \cite[pp.~481--482]{kronecker:01}
\end{quote}
Kronecker took the opportunity to clarify the methodological stance towards analysis that is appropriate to these issues. Even though, in his work, he avoided the notion of an ``arbitrary'' real number, he recognized the importance of understanding particular number systems in analytic terms. For example, given a symbolic representation of a real number, one may wish to compute rational approximations. But knowing that a nonnegative real number is nonzero is not the same as having a positive, rational lower bound on its values.\footnote{In modern terms, one can obtain such a lower bound by computing rational approximations until one obtains one that is sufficiently accurate to bound the number away from zero. Thus the statement ``if $r \neq 0$, then $|r| >0$'' was accepted by the Russian school of constructive mathematics in the 1950's and 1960's. This implication is equivalent to ``Markov's principle,'' which is, however, rejected by strict constructivists. See, for example, \cite{troelstra:van:dalen:88both}.} Kronecker mentioned the problem of bounding a nonzero determinant away from zero as an example of a problem that is generally difficult. He also noted that, given two convergent series, it can be difficult to determine which one has a greater value. These facts are now quite familiar in constructive and computable analysis.\footnote{See, for example, Troelstra and van Dalen \cite{troelstra:van:dalen:88both}.}

Kronecker's treatment of characters provides an interesting combination of Dirichlet's approach and modern ones. Fixing a modulus $m$, Kronecker, like Dirichlet, provided a fully explicit description of the characters modulo  $m$ in terms of primitive elements of the powers of primes giving $m$ and primitive roots of unity. However, with judicious choice of notation, he was at the same time able to suppress extraneous detail. For example, fixing primitive elements $\gamma, \gamma_0, \gamma_1, \ldots, \gamma_g$, one can express any element $r$ relatively prime to $m$ in the form
\[
 r \equiv \gamma^\rho \gamma_0^{\rho_0} \cdots \gamma_g^{\rho_g} \bmod{m}.
\]
But then one can ``package'' the representing values $(\rho, \rho_0, \rho_1, \ldots, \rho_g)$ as a tuple, the ``index system of $r$,'' denoted $\mathrm{Indd}\ r$. For each of the cyclic groups in the decomposition of the group of units modulo $m$, suppose we also choose corresponding roots of unity, $\omega, \omega_0, \omega_1, \ldots, \omega_g$.\footnote{As above, the powers of two require special treatment. In the discussion, Kronecker assumes that $m$ is divisible by 8, in which case the units modulo that power of two form a product of two cyclic groups; the values of $\rho$ and $\rho_0$ are the indices in those two groups, and $\omega$ and $\omega_0$ are the corresponding roots of unity.}
 \begin{quote}
 Now let $r$ be a unit modulo $m$, and $\mathrm{Indd}\ r = (\rho, \rho_0, \rho_1, \ldots, \rho_g)$; we now assign to $r$ the root of unity:
\[
 \Omega(r) = (-1)^\rho \omega_0^{\rho_0} \omega_1^{\rho_1} \ldots \omega_g^{\rho_g}
\]
 which we call a \emph{character} of $r$, since the index system $(\rho, \rho_0, \rho_1, \ldots, \rho_g)$, and hence $\Omega(r)$, are uniquely determined.\footnote{Es sie nun $r$ eine Einheit modulo $m$, und $\rm{Indd}\ r = (\rho, \rho_0, \rho_1, \cdots, \rho_g)$; ordnen wir $r$ jetzt die Einheitswurzel:
\[
 \Omega(r) = (-1)^\rho \omega_0^{\rho_0} \omega_1^{\rho_1} \ldots \omega_g^{\rho_g}
\]
zu, so geh\"ort zu jeder Einheit $r$ eine und nur eine Einheitswurzel $\Omega(r)$, welche wir einen \emph{Charakter} von $r$ nennen wollen, denn durch $r$ is ja das Indexsystem $(\rho, \rho_0, \cdots)$, also $\Omega(r)$ eindeutig bestimmt.} \cite[p.~444]{kronecker:01}
\end{quote}
We obtain \emph{all} the possible characters by fixing primitive roots of unity $\omega, \omega_0, \ldots, \omega_g$, so that every tuple of roots can be represented in the form
\[
 \omega^k, \omega_0^{k_0}, \omega_1^{k_1} \ldots, \omega_g^{k_g},
\]
 where the $k$ and $k_i$'s are less than the cardinality of the corresponding cyclic group.
\begin{quote}
When there is no fear of misunderstanding, we will, in the following, denote the underlying exponent system $(k, k_0, k_1, \cdots)$ by $(k)$ for short, and we will denote the corresponding character simply by
\[
\Omega^{(k)}(r).
\]
Here, again, to each system of values $(k, k_0, \cdots)$ and each unit $r$ there clearly corresponds a character $\Omega^{(k)}(r)$.\footnote{``Wenn kein Mi{\ss}verst\"andnis zu bef\"urchten ist, wollen wir im folgenden das zu Grunde gelegte Exponentensystem $(k, k_0, k_1, \ldots)$ kurz durch $(k)$ und den zugeh\"oringen Charakter einfacher durch
\[
\Omega^{(k)}(r)
\]
bezeichnen. Auch hier entspricht f\"ur ein festes Wertsystem $(k, k_0, \cdots)$ jeder Einheit offenbar $r$ ein Character $\Omega^{(k)}(r)$.''} \cite[p.~445]{kronecker:01}
\end{quote}
In contrast, say, to Hadamard, the index $(k)$ is not arbitrary; rather, $(k)$ is taken to range over the specific representing data. But the notation and organization of the proof enables us to ignore the details of the representation where they are not needed. Kronecker noted immediately that $\Omega^{(k)}(r)$ depends only on the value of $r$ modulo $m$, since any two values with the same residue have the same index. He also notes that we have $\Omega^{(k)}(rr') = \Omega^{(k)}(r) \Omega^{(k)}(r')$. The choice of representation has the further nice property that $\Omega^{(k)}(r) \Omega^{(k')}(r) = \Omega^{(k + k')}(r)$ where $(k + k')$ denotes the result of adding the elements of the tuples $(k)$ and $(k')$, modulo the cardinality of the associated cyclic groups. Kronecker presented the first orthogonality principle:
\[
 \sum_{(r)} \Omega^{(0)}(r) = \ph(m),
\]
where $r$ ranges over a system of residues of the units modulo $m$ and $(0)$ is the index of the trivial character, and
\[
  \sum_{(r)} \Omega^{(k)}(r) = 0
\]
for the remaining characters. Kronecker's proof requires unfolding the notation and calculating, but thereafter the fact can be recalled and used in the above form. Similarly, he expressed the second orthogonality relation by the equations
\[
 \sum_{(k)} \Omega^{(k)}(r_0) = \ph(m),
\]
when $r_0$ is congruent to $1$ modulo $m$, and
\[
 \sum_{(k)} \Omega^{(k)}(r) = 0
\]
otherwise. He also expressed the dependence of a Dirichlet series on a character in terms of a dependence on $(k)$, by writing
\[
 L^{(k)}(z) = \sum_{n = 1}^\infty \frac{\Omega^{(k)}(n)}{n^z} = \prod_p \frac{1}{1 - \frac{\Omega^{(k)}(p)}{p^z}}.
\]

Kronecker's presentation provides us with an important lesson. The exposition is clearly designed to bring the central ideas to the fore and highlight relevant information, but succeeds in doing so while providing explicit representations and algorithms throughout. This shows that although the moves towards abstraction that we have documented in this section can sometimes \emph{make it possible} to suppress or even ignore algorithmic information, they do not \emph{require} one to do so. In other words, the conceptual reorganization opens the door to the use of other means of describing mathematical objects and operations on them, means that can supplant algorithmic aspects. The question as to whether such means are permissible, meaningful, and appropriate to mathematics lay at the core of twentieth century foundational debates.

\subsection{Landau}
\label{landau:section}

Born in 1877, Edmund Landau received his doctorate at the University of Berlin in 1899, having studied number theory under Frobenius. He completed a \emph{Habilitation} thesis on Dirichlet series in 1901. Landau later presented proofs of Dirichlet's theorem in two textbooks, his 1909 \emph{Handbuch der Lehre von der Verteilung der Primzahlen} \cite{landau1909} and his 1927 \emph{Vorlesungen \"{u}ber Zahlentheorie} \cite{landau:27}.  We shall here focus on his presentation in the 1909 work, since the later presentation is already essentially what we have portrayed as ``contemporary'' in Section~\ref{contemporary:section}.

Landau began by introducing the characters via the construction in terms of primitive roots and roots of unity.  However, the notation that he used to denote them changed over the course of the book.  Initially, he used a notation that is quite similar to Dirichlet's notation for $L$-functions, writing
\[
\chi_{(a_{1}, a_{2}, \ldots, a_{r}, a, b)}(n),
\]
where the index system $a_{1}, a_{2}, \ldots, a_{r}, a, b$ serves to distinguish the characters just as Dirichlet's $\mathfrak{a}, \mathfrak{b}, \mathfrak{c}, \mathfrak{c'}, \ldots$ distinguished his $L$-functions.  However, after proving that there are $\ph(d)$ such characters, he simplified his notation to
\[
\chi_{1}(n), \chi_{2}(n), \ldots, \chi_{\ph(d)}(n),
\]
and to $\chi_{x}(n)$ in the general case, thus adopting notation similar to both de la Vall\'{e}e-Poussin and Hadamard.

But whereas Landau constructed the characters in the same way as Dirichlet, Dedekind, Hadamard and de la Vall\'{e}e-Poussin, he also recognized that they are characterized by certain key \emph{properties}, and that it is only these properties that are needed in the proof.  Indeed, after introducing the abbreviations for the characters, Landau wrote:
\begin{quote}
\ldots I will prove four short and elegantly-worded theorems about them [the characters]. The reader may then quickly forget the rather complicated definition of these functions completely, and need only remember that the existence of a system of $h$ distinct functions which possesses the four properties has been proved.\footnote{``\ldots ich werde \"{u}ber sie vier S\"{a}tze mit sehr kurzem und elegantem Wortlaut beweisen.  Alsdann darf der Leser bald die recht komplizierte Definition dieser Funktionen vollkommen vergessen und braucht sich nur zu merken, da{\ss} die Existenz eines Systems von $h$ verschiedenen Funktionen bewiesen worden ist, welche die vier Eigenschaften besitzen.''} \cite[p.~404]{landau1909}
\end{quote}

The four theorems that Landau was referring to are the following:
\begin{quote}
{\bf Theorem 1:} For any two positive numbers $n$ and $n'$, 
\[
\chi(nn')=\chi(n)\chi(n').
\]  
This ``law of multiplication'' holds for each of the $h$ functions \ldots

{\bf Theorem 2:} For $n \equiv n' \pmod{k}$,
\[
\chi(n)=\chi(n')
\]
\ldots

{\bf Theorem 3:} When $n$ runs through a complete residue system modulo $k$, for $x=1$, i.e.\ for the principal character $$\sum_{n}\chi_{x}(n)=h,$$ however for $x=2,\ldots h$, i.e.\ for all other characters 
\[
\sum_{n}\chi_{x}(n)=0
\]
\ldots

{\bf Theorem 4:} When $n$ is fixed and the sum
\[
\sum_{x=1}^{h}\chi_{x}(n)
\]
extends over all $h$ functions, then 
\[
\sum_{x=1}^{h}\chi_{x}(n)=h \ \mbox{for} \ n \equiv 1 \pmod{k},
\]
\[
\sum_{x=1}^{h}\chi_{x}(n)=0 \ \mbox{for} \ n \not\equiv 1 \pmod{k},
\]
therefore for all $k-1$ other residue classes modulo $k$.\footnote{``{\bf Satz 1:} Es ist f\"{u}r zwei ganze positive Zahlen $n, n'$ 
\[
\chi(nn')=\chi(n)\chi(n').
\] 
Von jeder der $h$ Funktionen wird also dies ``Multiplikationsgesetz'' behauptet \ldots

\noindent {\bf Satz 2:} Es ist f\"{u}r $n \equiv n' \pmod{k}$
\[
\chi(n)=\chi(n')
\]
\ldots

\noindent {\bf Satz 3:} Wenn $n$ ein vollst\"{a}ndiges Restsystem modulo $k$ durchl\"{a}uft, ist f\"{u}r $x=1$, d.h. f\"{u}r den Hauptcharakter $$\sum_{n}\chi_{x}(n)=h,$$ dagegen f\"{u}r $x=2, \ldots, h$, d.h. f\"{u}r alle \"{u}brigen Charaktere 
\[
\sum_{n}\chi_{x}(n)=0
\]
\ldots

\noindent {\bf Satz 4:} Wenn $n$ festgehalten und die Summe 
\[
\sum_{x=1}^{h}\chi_{x}(n)
\] 
\"{u}ber alle $h$ Funktionen erstreckt wird, so ist $$\sum_{x=1}^{h}\chi_{x}(n)=h \ \mbox{f\"{u}r} \ n \equiv 1 \pmod{k},$$ dagegen $$\sum_{x=1}^{h}\chi_{x}(n)=0 \ \mbox{f\"{u}r} \ n \not\equiv 1 \pmod{k},$$ also f\"{u}r alle $k-1$ \"{u}brigen Restklassen modulo $k$.''} \cite[pp.~401--408]{landau1909}
\end{quote}
These are exactly the four ``important relations'' given by de la Vall\'ee-Poussin, and the analogues of the four properties enumerated by Dirichlet in 1840. What is novel here is Landau's explicit recognition that only these properties are used in the proof, and that the specific construction is only needed to show the existence of a system of functions that satisfy them.

Like Hadamard, Landau did not go so far as to allow the characters themselves to index the sum in theorem Theorem 4 above, relying  on a natural number proxy. But, like de la Vall\'ee-Poussin, he gave an extensional classification of the three types of characters, though he also gave an additional description of the real characters in intensional terms.

Just as Landau's notation for characters changed over the course of the book, so, too, did his notation for $L$-functions. Initially, he denoted them in a manner similar to Hadamard's, writing $L_{x}(s)=\sum_{n=1}^{\infty}\frac{\chi_{x}(n)}{n^{s}}$.  However, later on, after concluding his presentation of the proof of Dirichlet's theorem, Landau adopted a ``more convenient'' notation:
\begin{quote}
Now let 
\[
L_{x}(s) = \sum_{n=1}^{\infty}\frac{\chi_{x}(n)}{n^{s}}
\]
be the function corresponding to the character $\chi(n)=\chi_{x}(n)$; it is now more convenient to include the character in the notation, 
\[
L(s, \chi),
\] 
and, only when there is no fear of misunderstanding, write 
\[
L(s)
\]
for short, as before.\footnote{``Es sei nun 
\[
L_{x}(s) = \sum_{n=1}^{\infty}\frac{\chi_{x}(n)}{n^{s}}
\] 
die dem Charakter $\chi(n)=\chi_{x}(n)$ entsprechende Funktion; es ist jetzt bequemer, um die Charakter in die Bezeichnung aufzunehmen, 
\[
L(s, \chi)
\] 
zu schreiben, und nur, wenn kein Mi{\ss}verst\"{a}ndnis zu bef\"{u}rchten ist, wie fr\"{u}her kurz 
\[
L(s).\mbox{''}
\]} \cite[p.~482]{landau1909}
\end{quote}
Landau did not explain \emph{why} the modern notation is more convenient, but we will offer some suggestions in the next section.

Landau's 1909 presentation is interesting because it has a transitional feel. Although he gave an intensional construction of the characters, he not only provided a thorough enumeration of their properties, but went out of his way to emphasize that these properties are all that matters to the proof. This is borne out by the fact that the division of characters into three classes can be carried out extensionally, without reference to their construction. And although he initially introduced natural number indices for the characters and $L$-functions and summed over these indices, he eventually adopted the modern functional notation for $L$-functions.

When he wrote his 1927 textbook \cite{landau:27}, Landau finally made the transition to a proof that is extremely close to our contemporary version.  Indeed, he no longer defined the characters by describing how they are constructed, but, rather, defined them in terms of their characteristic properties.  Moreover, he summed over the characters themselves, using expressions such as $\sum_{\chi}\chi(a)$. The only real difference between his proof and the one we presented in Section \ref{contemporary:proofs:section} is that he did not develop the general notion of a group-theoretic character, but, rather, defined them in terms of the \emph{particular} group of units modulo $m$. This makes sense, given that the work is an elementary number theory textbook and characters are not needed for any other purpose. Roughly speaking, if we combine Landau's 1927 presentation with Weber's 1882 more general treatment of characters, the result is the presentation in Section~\ref{contemporary:section}.

\section{Changes in mathematical method}
\label{analysis:section}

In Section~\ref{objects:section}, we considered a number of ways in which characters are treated as bona-fide objects in contemporary proofs of Dirichlet's theorem. We noted that none of these features are present in Dirichlet's original proof, for the simple reason that Dirichlet did not isolate or identify the notion of a character. Dedekind's 1863 presentation of Dirichlet's theorem did not use the term ``character,'' but he did introduce the notation, $\chi(n)$, for characters, and identified the defining property of a homomorphism as their ``characteristic property.'' As noted in Section~\ref{dedekind:section}, after Weber's paper of 1882 the general notion of a group character was in place, and all the subsequently published proofs of Dirichlet's theorem use the terminology of characters.

But now the set of characters modulo $m$ can be defined extensionally, as the set of nonzero homomorphisms from $(\mathbb Z/ m \mathbb Z)^*$ to the complex numbers, or intensionally, as functions defined by certain algebraic expressions involving certain primitive elements modulo the prime powers occurring in the factorization of $m$, and certain complex roots of unity. Even though the two definitions give rise to the same set of characters, proofs can differ in the extent to which they rely on the specific representations or the abstract characterizing property. Dirichlet's proof relied only on the symbolic representations. Somewhat surprisingly, both Dedekind's and Hadamard's division of the characters into the trivial, real, and complex cases was also described in terms of the characters' representations, even though the distinction is naturally expressed in terms of the values they take. Kronecker and de la Vall\'ee-Poussin provided both descriptions, and even though Kronecker made it clear that all operations and classifications can be carried out, algorithmically, in terms of the canonical representations, his careful choice of notation and organization made the extensional properties salient. By 1927, Landau clearly favored the extensional characterization in his textbook.

We have also observed that modern notation like $\sum_\chi \chi(n)$ allows us to carry out summations over the finite set of characters modulo $m$, but that this notation was not used by Dirichlet's early expositors. Dirichlet, Dedekind, and Kronecker all took summations to range over the tuples of integers representing the characters via an explicit algebraic definition, though Kronecker's way of letting a variable $(r)$ range over these tuples is more attractive than Dirichlet's use of $\mathfrak{a}, \mathfrak{b}, \mathfrak{c}, \mathfrak{c'}, \ldots$. In 1882, Weber wrote his sums with ellipses. Curiously, Hadamard in 1896, and Landau in 1909, assigned arbitrary integer indices to the characters, and took sums to range over those indices. Of the proofs of Dirichlet's Theorem~\ref{dirichlet:theorem} that we have considered, the only one that takes summation over characters at face value is that of de la Vall\'ee-Poussin, who nonetheless adopted Dirichlet's 1840 notation $S_\chi$ to denote such sums. By 1927, however, Landau had adopted the modern notation.

Finally, we have emphasized the modern tendency to represent the dependence of an $L$-series on a character $\chi$ as a functional dependence, with the notation $L(s,\chi)$. Once again, de la Vall\'ee-Poussin was the only nineteenth century expositor of Dirichlet's theorem to do so, with the notation $Z(s,\chi)$. We saw that Landau made the transition in the middle of his 1909 book, and that in his 1927 textbook he relied exclusively on the notation $L(s,\chi)$.

At issue in all these developments is whether characters could be treated in much the same way as natural numbers and real numbers, or whether characters are different sorts of objects, whose treatment has to be mediated by more ``concrete'' mathematical representations. We contemporary readers of the nineteenth century literature are now so familiar with a modern perspective that it can be hard for us to appreciate the reasons it took so long for the community to adopt it. Let us begin our analysis of the history, then, by reflecting on the countervailing pressures.

We have had a lot to say about the notational and expressive conventions at play. The difficulty of settling on stable and useful conventions should not be minimized. For example, any mathematician writing a proof has to choose names for variables: should one use $m$, $n$, and $k$ to range over natural numbers, or $x$, $y$, and $z$? One desiderata is to maintain continuity with the background literature, but other constraints come into play: for example, one reason to favor $m$, $n$, and $k$ may be that $x$, $y$, and $z$ are natural choices to range over other objects, like real or complex numbers, arising in the proof. It is significant, though not surprising, that Dedekind, 26 years after Dirichlet's proof was published, and de la Vall\'ee-Poussin, fully 60 years later, both stuck with Dirichlet's choice of $m$ and $k$ for the first term and common difference in the arithmetic progression, respectively. Even today, $q$ is often used to range over prime numbers in the definition of the $L$-series, and we still use the letter $L$ exclusively, though most contemporary proofs use $a$ and $d$ in place of $m$ and $k$.

In a similar way, one has to settle on choices of notation. This also requires thought, even in places where the norms that govern the notation are fairly clear. For example, assuming one divides the set $C$ of characters modulo $k$ into the set consisting of the trivial character, $C_{\mathit{triv}}$, the set of real characters, $C_{\mathit{real}}$, and the set of complex characters, $C_{\mathit{complex}}$, it is clear that a sum over $C$ can be broken up accordingly:
\[
\sum_{\chi \in C} F(\chi) = \sum_{\chi \in C_{\mathit{triv}}} F(\chi) +
  \sum_{\chi \in C_{\mathit{real}}} F(\chi) +
  \sum_{\chi \in C_{\mathit{complex}}} F(\chi).
\]
But, even so, one has to settle on the notation to express this relationship, and, as the history of Dirichlet's theorem shows, it can be a long time before a particular means of expression becomes standard. The modern notation emphasizes the parallels between summing over sets $C$ of characters and summing over sets $C$ of natural numbers, but it is by no means obvious that conflating the two is a good idea.

In Section~\ref{generalization:section}, we briefly touched on Frege's treatment of functions in the logical system of his \emph{Begriffsschrift} \cite{Frege1879}, developed in greater detail in \emph{Grundgesetze der Arithmetik} \cite{frege:grundgesetze}. Contemporary readers of the \emph{Grundgesetze} may be struck by how many pages are devoted to explaining the syntax of the formal language. For those familiar with modern logic, Frege's lengthy explanations seem fiddly and pedantic. But, for Frege, getting the grammatical rules worked out was much of the battle. The syntax of a language only seems trivial when you already know how to speak it.

But it would be a mistake to suggest that all the considerations that Dirichlet's successors faced were ``merely'' notational. Means of expression often rely on substantial features of our understanding of the nature of that which is expressed. For example, consider the modern way of writing Dirichlet's famous example of 1927:
\[
f(x) = \begin{cases}
  a & \mbox{if $x$ is rational} \\
  b & \mbox{otherwise.}
\end{cases}
\]
This bears superficial similarity to case-based definitions of number-theoretic functions, or piecewise definitions of analytic functions, that were familiar at the time. But it was novel to base a case distinction on the property of being rational, and, indeed, the notation masks significant assumptions about what it means to define a function on the reals. Consider Frege's brief account of the evolution of the function concept in 1891:
\begin{quote}
In the first place, the field of mathematical operations that serve for constructing functions has been extended.  Besides addition, multiplication, exponentiation, and their converses,  the various means of transition to the limit have been introduced---to be sure, people have not always been clearly aware that they were thus adopting something essentially new.  People have gone further still, and have actually been obliged to resort to ordinary language, because the symbolic language of Analysis failed; e.g. when they were speaking of a function whose value is 1 for rational and 0 for irrational arguments.\footnote{``Erstens n\"{a}mlich ist der Kreis der Rechnungsarten erweitert worden, die zur Bildung einer Funtion beitragen.  Zu der Addition, Multiplikation, Potenzierung und deren Umkehrungen sind die verschiedenen Arten des Grenz\"{u}berganges hinzugekommen, ohne da{\ss} man allerdings immer ein klares Bewu{\ss}tsein von dem wesentlich Neuen hatte, das damit aufgenommen werde.   Man ist weiter gegangen und sogar gen\"{o}tigt worden, zu der Wortsprache seine Zuflucht zu nehmen, da die Zeichensprache der Analysis versagte, wenn z.B. von einer Funktion die Rede war, deren Wert f\"{u}r rationale Argumente 1, f\"{u}r irrationale 0 ist.''} \cite[p.~12]{Frege91}
\end{quote}
Frege was not merely concerned to have a convenient notation to express the definition by cases. The definition of $f$ above does that perfectly well, even today. The point is rather that the notation should come with clear rules of use. That is what Frege took to be lacking in the casual use of ordinary language, and what he took his formal system to provide. Similar methodological concerns lie beneath the surface whenever we write $\lim_{n \to \infty} a_n$ to denote the limit of a sequence of real numbers, or $I + J$ to denote the sum of two ideals in a ring of algebraic integers. There is nothing tame about the infinitary operations that underlie the notation. 

In contrast, foundational questions regarding the use of characters may seem mild. After all, it is easy to represent the characters modulo a positive integer $k$, and any character is determined by the values it takes on the finitely many residues modulo $k$. Nonetheless, the use and treatment of characters in proofs of Dirichlet's theorem bears upon central questions regarding the use and treatment of functions more generally, specifically regarding the relationship between a function and its various representations.

Let us think of Dirichlet $L$-functions $L(s,\chi)$ in quasi-computational terms. Such a function should take, as input, a real (or complex) value $s$ and a character $\chi$, and return a complex number. Let us set aside the (important) question as to what it means to take a real or complex number as input, or return one as output. What does it mean to accept a character as an input? Should one think of the character as being ``presented'' to the function as an infinite set of input-output pairs? Or as the list of finite values on the residues modulo $k$? (In that case, does the same conception work for function arguments that are not finitely determined?) Should one think, rather, of $\chi$ as being some sort of procedure, or subroutine? If so, what sorts of procedures and subroutines are allowed in the definition of a functional with argument $\chi$? Perhaps, instead, we should identify $\chi$, as Dirichlet did, with its representation in terms of an expression involving certain roots of unity. But recall that those representations relied on choices of primitive elements in the representation of $k$, although it turns out that the value of the $L$-function does not vary with different choices of the primitive elements. On our conception, does $L$ somehow ``depend'' on the choices of primitive elements?

These issues arise not only with respect to functional notation, but also with respect to statements involving quantifiers over characters. For example, the proof of the second orthogonality lemma requires the following fact:
\begin{quote}
If $n$ is relatively prime to $k$, there is a character $\chi$ such that $\chi(n) \neq 1$.
\end{quote}
If we prove this as a separate lemma, we can then invoke the lemma for a given $k$ to ``obtain'' a $\chi$ with the relevant property. But what exactly have we ``obtained''? A table of values? A procedure? A representation? If the lemma constructs a $\chi$ via a representation of one sort, but our proof of the second orthogonality lemma relies on a different representation of characters, is it legitimate to apply the lemma in that context?

The modern theory of computability and the semantics of programming languages offers various ways of thinking about computer programs which take functions are arguments. The issues are subtle and complex. In contrast, modern mathematics followed a route whereby these subtleties are for the most part set aside. In particular, they are deemed incidental to the proof of Dirichlet's theorem. Roughly speaking, to make sense of a functional $F(f)$ with function argument $f$, set theory identifies $f$ ``canonically'' with its extension, a set of input-output pairs, without concern as to how $f$ is represented or how (and whether) one can ``compute'' $F(f)$. Set theory then imposes on mathematical language the restriction that the definition of such a functional $F$ can only depend on the extension of $f$. (Or, put differently: modern mathematical conventions evolved to ensure the latter fact, and axiomatic set theory was designed to model and explain those conventions.)

In the nineteenth century, the answers to the questions raised above were not at all obvious. Indeed, they are still debated among logicians and foundationally-minded mathematicians today. Even for those inclined to dismiss those questions as irrelevant to the proof of Dirichlet's theorem, it would not have been immediately clear as to whether they really could be dismissed, and, if so, how that should be done. What may seem to be ``merely notational'' developments in the presentation of Dirichlet's theorem were part and parcel of the broader mathematical community's attempt to fashion an understanding of functions as objects that would better support the mathematics of the time.

Setting aside issues of meaning, there may also be concerns about \emph{correctness}. Any choice of notation that draws on an analogy between different domains presupposes that the analogy is appropriate, which is to say, that sufficiently many properties carry over, and that there are sufficient safeguards to bar the ones that do not. For example, a sum of the form $\sum_{x \in S} t(x)$ generally makes sense when $S$ is finite and the term $t(x)$ takes values on a domain with an addition that is associative and commutative, so that the order the one ``runs through'' the elements of $S$ does not matter. Summing over infinite totalities is more subtle; it typically depends on having a notion of convergence for the domain in which $t(x)$ takes values, and worrying about the order in which terms are summed. Now, choosing the notation $\sum_\chi t(\chi)$ may not be a good idea if it encourages users to transfer facts and properties between domains in an invalid way. This would explain why Dirichlet and de la Vall\'ee-Poussin chose a special notation $S_\chi$ for such sums, as such a new symbol would not come with any unwanted baggage. Adopting the contemporary notation may have required some kind of assurance that the intended contexts would be sufficient to control for proper use.

In sum, mathematical conventions evolve, expand, and change. Any time a mathematician writes a line of text, he or she is situated in a tradition with implicit norms and conventions, and the line of text just written becomes part of that tradition. The fact that so many of these conventions are communicated implicitly does not make them any less important to mathematical understanding.

When a mathematician writes a proof, the intent is that the inferences contained therein will be deemed by his or her colleagues to be correct and justified. Where the inferences are instances of familiar patterns, one desires that they will be easily recognized as such, so that the reader's effort is conserved for more substantial cognitive tasks. When the inferences rely on assumptions that may be considered dubious, or push familiar patterns of reasoning into unexplored territory, there is greater concern not only as to whether the reasoning will be recognized as correct, but also as to whether it will be deemed \emph{appropriate} to addressing the mathematical issues at hand. Thus there are always strong pragmatic pressures to stick close to established convention, and one should not expect fundamental aspects of the language and methods of mathematics to change in novel ways, unless there are strong forces pushing for such change.

In the case of Dirichlet's theorem, we have seen some of the ways in which Dirichlet's successors chose to modify, or ``improve,'' his presentation. Now let us try to understand some of the perceived benefits. Kronecker's proof was explicitly designed to fill in information that was absent from Dirichlet's proof, in the form of explicit bounds on a quantity asserted to exist. Many of the other developments were explicitly designed towards paving the way to useful generalizations. This was clearly Dedekind's intent, for example, in pointing out that the Euler product formula depends only on certain multiplicative properties of the terms occurring in the sum. Similarly, abstracting proofs of properties of characters on $(\ZZ / m\ZZ)^*$ that rely on features specific to the integers modulo $m$ paves the way to extending these properties to group characters more generally. Characters and their properties form the basis for representation theory, which has been an essential part of group theory since the turn of the twentieth century \cite{hawkins:71,mackey:80}. Authors like Hadamard, de la Vall\'ee-Poussin, and Landau were also interested in extending Dirichlet's methods to other kinds of Dirichlet series, which now play a core role in analytic number theory.

But it would be a mistake to attribute all the benefits of the expository innovations we have considered to increased generality. After all, these innovations play an equally important role in fostering a better understanding of Dirichlet's proof itself, by highlighting key features of the concepts and objects in question, motivating the steps of the proof, and reducing cognitive burden on the reader by minimizing the amount of information that needs to be kept in mind at each step along the way. It is true that these benefits often support generalization, but they do so in part by making our thinking \emph{vis-\`a-vis} Dirichlet's proof itself more efficient.

For example, proofs of Dirichlet's theorem that rely on explicit representations of the characters require us to keep the details of the representation in mind. Recall Dirichlet's original representation of an arbitrary character:
\[
 \chi(n) = \theta^{\alpha} \ph^{\beta} \omega^{\gamma} \omega'^{\gamma'} \cdots
\]
Here the reader has to keep in mind that $\alpha$, $\beta$, $\gamma$, $\gamma'$, and so on are the indices of $n$ with respect to primitive elements chosen in the decomposition of the abelian group $(\mathbb{Z}/m \mathbb{Z})^*$, and $\theta$, $\ph$, $\omega$, $\omega'$, and so on are corresponding roots of unity.
This information has to be kept in mind throughout the proof, because the nature of the objects and the dependences of $\alpha, \beta, \gamma, \gamma'$ on $n$ may play a role in licensing an inference or calculation.
Moreover, recall that we obtain all the characters by expressing all the roots
\[
\theta = \Theta^\mathfrak{a}, \ph = \Phi^\mathfrak{b}, \omega = \Omega^\mathfrak{c}, \omega' = \Omega^{\mathfrak{c}'}, \ldots,
\]
in terms of primitive roots of unity $\Theta$, $\Phi$, $\Omega$, $\Omega'$, \ldots, and letting $\mathfrak{a}, \mathfrak{b}, \ldots$ range over the appropriate exponents. Once again, this information has to be remembered throughout. Later proofs are easier to read simply because they do not require us to keep as much information in mind, and highlight the relevant dependences when they are needed.

One way of achieving this is by reorganizing the proof in such a way that some of the relevant information is localized to particular facts and calculations. For example, even if one resorts to representations to prove the orthogonality relations, if this is the only place they are used, then they do not need to be ready to hand when these relations are invoked in a calculation later on, where other analytic expressions and their properties are the objects of focus. Thus modularity reduces cognitive burden, and makes it easier to keep track of the global structure of the argument, providing high-level outlines, or sketches, of the proof. Such restructuring paves the way to generality: isolating context-specific details in well-insulated modules means that one can adapt the proof by changing the modules while preserving their external interfaces. But, to repeat, generality is not the only benefit: the restructuring improves the readability of the original proof as well.

Even in situations where there is a lot of information in play at once, judicious notation and means of expression can make important relationships between the data more salient. For example, in Section~\ref{dirichlet:section}, we saw that Dirichlet wrote
\begin{multline*}
\sum\frac{1}{q^{1 + \rho}} + \frac{1}{2}\sum\frac{1}{q^{2 +2 \rho}} + \frac{1}{3}\sum\frac{1}{q^{3 + 3 \rho}} + \ldots \notag \\
= \frac{1}{K}\sum \Theta^{-\alpha_{m}\mathfrak{a}}\ \Phi^{-\beta_{m}\mathfrak{b}}\Omega^{-\gamma_{m}\mathfrak{c}}\Omega^{-\gamma^{'}_{m'}\mathfrak{c}'} \ldots \log L_{\mathfrak{a},\mathfrak{b},\mathfrak{c},\mathfrak{c}'\ldots.}
\end{multline*}
where we would write
\[
\sum_{\chi \in \widehat{(\mathbb{Z} / q \mathbb{Z})^*}}\overline{\chi(m)}\log L(s, \chi)=\ph(q) \sum_{p \equiv m \pmod{q}}\frac{1}{q^{s}} \ + \ O(1).
\]
The second presentation highlights the relationship between the series $L(s, \chi)$ and the characters; in particular, the value $L(s,\chi)$ is multiplied by the value of the conjugate character, $\overline \chi$, at $m$ in each term of the sum on the left-hand side. Dirichlet's notation obscures this relationship, since the logarithm of the appropriate $L$-series $L_{\mathfrak{a}, \mathfrak{b}, \mathfrak{c}, \mathfrak{c}'\ldots}$ is multiplied by the character $\Theta^{-\alpha_{m}\mathfrak{a}}\ \Phi^{-\beta_{m}\mathfrak{b}}\Omega^{-\gamma_{m}\mathfrak{c}}\Omega^{-\gamma^{'}_{m'}\mathfrak{c}'}\ldots$ and the reader must sift through this expression to notice that the exponents $-\mathfrak{a}, -\mathfrak{b}, -\mathfrak{c}, -\mathfrak{c}', \ldots$ correspond to the subscripts of the $L$-series.

The same problem applies to Hadamard's method of using arbitrary indices for the characters: given a character $\psi_v$, there is no natural name for the index of its conjugate $\overline{\psi}_{v}$. Of course, one can introduce such a notation, but that requires keeping track of the relationship between the two notations, that is, the conjugation of characters and the associated operation on indices. And grouping the characters into different classes requires grouping the indices into different classes, again yielding an uncomfortable duality. (At least Kronecker maintained a monist consistency, insisting that all operations on characters are operations on their representations. So when $\Omega^{(k)}$ is a character, Kronecker could write $\Omega^{(-k)}$ for its conjugate, since the conjugate character is obtained by negating the elements of the corresponding tuple.)

Similar considerations may explain why Landau changed his notation for $L$-functions in the middle of his 1909 work, from $L_{x}(s)$ to $L(s, \chi)$.  In the paragraph following his notation change, Landau showed that the theory of $L$-functions can be reduced to $L$-functions that correspond to particular types of characters, called \emph{proper characters}.  Roughly, proper characters modulo $k$ are those which cannot be obtained as a character modulo $K$ where $K<k$.  To show that the theory of $L$-series can be reduced in the appropriate way, he proved that if $\chi$ is an improper character modulo $k$ and $X$ is a corresponding proper character modulo $K$, we have \cite[482-483]{landau1909}
\[
L(s, \chi) = \prod_{\nu=1}^{c}\left(1-\frac{\varepsilon_{\nu}}{p^{s}_{\nu}}L_{0}(s, X)\right).
\]
Here the $\varepsilon_{\nu}$ are certain roots of unity, $c$ is a natural number (which can be 0, depending on how many prime factors of $k$ are contained in $K$), and the real part of $s$ is assumed to be greater than $1$. If Landau had kept his original notation, the left hand side of the above equation would be written as $L_{x}(s)$ where $\chi=\chi_{x}$.  But how would we represent the right hand side and, in particular, what index should we choose for $X$?  We should note that one obvious choice for an index, $x'$, would not be available, since Landau used this previously in connection with the conjugates of characters. In the text, Landau had identified a relationship between a character, $\chi$, and the corresponding proper character, $X$; having to translate this to a relationship between indices would only clutter the exposition. These issues are compounded when, in later sections, Landau wanted to obtain functional equations to relate $L(s, \chi)$ and $L(1-s, \overline{\chi})$ and in doing so referred to the distinction between proper and improper characters (see e.g. \cite[$\S$130]{landau1909}).  Using the old notation, we would need to keep track of three different subscripts to index the various characters.

Consider, finally, the issue of uniformity of notation. We have already discussed concerns associated with the notation $\sum_{\chi}$ for summation over characters. But there is an obvious benefit to using the same notation for summing over finite sets of characters and summing over finite sets of integers, namely, that the two operations really do share common properties. Indeed, the transfer is immediate via the Hadamard trick of assigning an integer index to each character. The option of introducing a new symbol every time one needs to index sums by finite sets of objects is clearly untenable, as it would result in a confusing explosion of notations. It does not seem at all surprising that Dirichlet's 1841 notation $S_\chi$ was short-lived.

To sum up, then, we have identified a number of advantages to the rewritings of Dirichlet's proof considered in the last section:
\begin{itemize}
 \item Essential properties of key objects (or expressions) were isolated, reducing the amount of information that someone reading the proof has to keep in mind at each step.
 \item The organization of proofs became more modular, with information localized to very specific parts of the proof.
 \item Expressions became more readable, since irrelevant details were suppressed, and the features that remained made dependences and relationships between terms more salient.
 \item Notation became more uniform, highlighting commonalities between different domains.
\end{itemize}
Changes like this often go hand in hand with attempts to generalize concepts and methods to other domains, since managing and controlling the volume of domain-specific detail tends to bring to the fore aspects of the proof that transcend these specifics. But they also contribute to a better understanding of Dirichlet's proof itself, and make the proof easier to read and reproduce from memory.

Benefits such as these are often dismissed as merely ``pragmatic'' or ``cognitive,'' but this downplays the fact that such considerations effectively shape and justify the norms that guide our mathematical practice. The history and philosophy of mathematics need to take them seriously.

\section{Conclusions}
\label{conclusions:section}

Our history of Dirichlet's theorem has highlighted one mathematical development that encouraged a higher-order treatment of functions, but, of course, there were others. One important example is the development of functional analysis. In 1859, George Boole published \emph{A Treatise on Differential Equations} \cite{boole:59}, in which he introduced the subject as a study of ``variable quantities'' subject to ``known'' relations between their differential coefficients. That work is noteworthy for its observation that differential operators can be viewed as algebraic expressions subject to certain laws. Integration was also viewed abstractly; in the expression $y = \int \ph(x) dx + c$,
\begin{quote}
the symbol $\int$ denotes a certain process of integration, the study of the various forms and conditions of which is, in a peculiar sense, the object of this part of the Integral Calculus. \cite[p.~2]{boole:59}
\end{quote}
Modern functional analysis takes the view that operations on functions like differentiation and integration can themselves be viewed as functions, defined over spaces of other functions. In 1887, Vito Volterra published a seminal paper, ``Sopra le funzioni che dipendono de altre funzioni'' (``On functions that depend on other functions'') \cite{volterra:87}, which adopted such a viewpoint. In 1901, Hadamard published \emph{Le\c{c}ons sur le calcul des variations}, which helped establish the foundations of functional analysis, and, in fact, introduced the term ``functional.'' In 1930, Volterra published an English translation of a series of lectures he had given at the University of Madrid in 1925 \cite{volterra:30}, which began with a lengthy discussion of the notion. After presenting particular examples, he said:
\begin{quote}
 We shall therefore say that a quantity $z$ is a \emph{functional of the function $x(t)$ in the interval $(a,b)$ when it depends on all the values taken by $x(t)$ when $t$ varies in the interval $(a,b)$}; or, alternatively, \emph{when a law is given by which to every function $x(t)$ defined within $(a,b)$ (the independent variable within a certain function field) there can be made to correspond one and only one quantity $z$, perfectly determined, and we shall write
\[
z = F \left|\left[ x  \stackrel[a]{b}{(t)} \right]\right|.                                                                                                                                                                                                                                                                                                                                                                                                                                                                                                                                                                                                                                                                                                                                                                                                                                                                                                                              \]
}

This definition of a \emph{functional} recalls especially the ordinary general definition of a \emph{function} given by Dirichlet. \cite[p.~4]{volterra:30}
\end{quote}
The chapter as a whole describes an outlook that is essentially the contemporary mathematical viewpoint, yet in a way that makes clear that the viewpoint was one that the mathematical community was still getting used to.\footnote{See also the Griffith Evans' helpful introduction to the 1959 Dover reprinting of Volterra's lectures \cite{volterra:30}.}

Twentieth-century foundational efforts served moreover to unify the function concept and provide a comprehensive framework to facilitate such higher-order treatment. In 1905, Borel, Baire, Lebesgue, and Hadamard engaged in their famous debate as to whether it makes sense to consider ``arbitrary functions,'' not given by any rule or law.\footnote{The ``five letters'' are translated as an appendix to Moore \cite{moore:82}, reproduced in Ewald \cite{ewald:96}, volume 2, pages 1077--1086.} In 1914, Felix Hausdorff published his \emph{Grundz\"uge der Mengenlehre} \cite{hausdorff:14}, dedicated to Georg Cantor. The book established a modern set-theoretic foundation, and used it to support the development of point-set topology. After discussing the notion of an ordered pair, Hausdorff gave the concept of function its modern definition:
\begin{quote}
\ldots we consider a set $P$ of such pairs, satisfying the condition that each element of $A$ occurs as the first element of one and only one pair $p$ of $P$. In this way, each element $a$ of $A$ determines a unique element $b$, namely, the one to which it is connected in a pair $p = (a, b)$. We denote this element associated to $a$, which is determined by and dependent on $a$, by
\[
 b = f(a),
\]
and we say that on $A$ (i.e.~for all elements of $A$) a \emph{unique function} of $a$ is defined. We view two such functions $f(a), f'(a)$ as equal when and only when the corresponding sets of pairs $P, P'$ are equal, that is, when, for \emph{each} $a$, $f(a) = f'(a)$.\footnote{``\ldots betrachten wir eine Menge $P$ solcher Paare, und zwar von der Beschaffenheit, da{\ss} jedes Element $a$ von $A$ in einem und nur einem Paare $p$ von $P$ als erstes Element auftritt. Jedes Element $a$ bestimmt auf diese Weise ein und nur ein Element $b$, n\"amlich dasjenige, mit dem es zu einem Paare $p = (a,b)$ verbunden auftritt; dieses durch $a$ bestimmte, von $a$ abh\"angige, dem $a$ zugeordnete Element bezeichnen wir mit
\[
 b = f(a)
\]
und sagen, da{\ss} hiermit \emph{in $A$} (d.\ h.\ f\"ur alle Elemente von $A$) eine \emph{eindeutige Funktion} von $a$ definiert sei. Zwei solche Funktionen $f(a)$, $f'(a)$ sehen wir dann und nur dann als gleich an, wenn die zugeh\"oringen Paarmengen $P, P'$ gleich sind, wenn also, f\"ur \emph{jedes} $a$, $f(a) = f'(a)$ ist.''} \cite[p.~33]{hausdorff:14}
\end{quote}

Thus, in the early decades of the twentieth century, the modern view of a function took root. A function could map elements of any domain to any other; one could specify a function by specifying any determinate law; functions could serve as arguments to other functions; and functions could serve as elements of algebraic structures and geometric spaces. 

Our study of the treatment of characters in number theory has focused on only one small part of the grand historical development that resulted in this way of thinking. Despite its narrow focus, the case study has illuminated some important factors that contributed to the modern view. We have emphasized the difficulties inherent in developing a coherent, rigorous way of treating functions as objects, as well as concerns as to how one could do so in a way that is appropriate to mathematics and preserves all the information that is essential to a rigorous argument. We have also explored some of the benefits associated with the methodological changes. Many of these accrue to the ability of the new language and notation to highlight central features and relationships of the objects in question, while suppressing details that make it harder to discern the high-level structure of a mathematical argument. We also emphasized the importance of representing mathematical concepts and objects in order to capture those features that are uniform across different domains of argumentation, so that these uniformities can be packaged and used in a modular way.

Finally, we have argued that understanding and evaluating the considerations that shape mathematical language and inferential practice is an important part of the history and philosophy of mathematics. Over time, mathematics evolves in such a way as to support the pursuit of distinctly mathematical goals. One of these is the goal of maintaining an inferential practice with clear rules and norms, one that allows its practitioners to carry out, communicate, and evaluate arguments that can become exceedingly long and complex. Another is the goal of promoting efficiency of thought, leveraging whatever features we can to extend our cognitive reach and transcend our cognitive limitations. The challenge remains that of developing appropriate ways of coming to terms with such ``pragmatic'' and ``cognitive'' considerations.\footnote{See Avigad \cite{avigad:06b,avigad:08b,avigad:10d} and the 
essays in Mancosu \cite{mancosu:08} for initial philosophical attempts in that direction, and Harel et al.~\cite{harel:et:al:08,harel:kaput:91} for a pedagogical perspective.}



\end{document}